\documentclass[11pt,reqno]{amsart}
\usepackage{hyperref}
\usepackage{color}
\usepackage{graphicx}
\usepackage{microtype}
\usepackage[alphabetic]{amsrefs}

\newtheorem{thm}{Theorem}[section]
\newtheorem{lem}[thm]{Lemma}
\newtheorem{prop}[thm]{Proposition}
\newtheorem{cor}[thm]{Corollary}

\newcommand{\newtheorembox}[2]{\newtheorem{x#1}[thm]{#2}
    \newenvironment{#1}{\begin{x#1}\pushQED{\qed}}{\popQED\end{x#1}}}

\theoremstyle{definition}
\newtheorem{defn}[thm]{Definition}
\newtheorembox{ex}{Example}
\newtheorembox{rmk}{Remark}
\newtheorembox{const}{Construction}

\newcommand{\arxiv}[1]{\href{http://arxiv.org/abs/#1}{arXiv:#1}}

\newcommand{\cA}{\mathcal A}
\newcommand{\CC}{\mathbb C}
\newcommand{\cD}{\mathcal D}
\newcommand{\cF}{\mathcal F}
\newcommand{\cG}{\mathcal G}
\newcommand{\cP}{\mathcal P}
\newcommand{\defi}[1]{\emph{#1}}
\renewcommand{\div}{\operatorname{div}}

\newcommand{\isom}{\cong}
\newcommand{\Jac}{\operatorname{Jac}}
\newcommand{\link}{\operatorname{link}}

\newcommand{\mult}{\operatorname{mult}}
\newcommand{\NS}{\operatorname{NS}}
\newcommand{\opp}{\operatorname{opp}}
\newcommand{\Pic}{\operatorname{Pic}}
\newcommand{\PP}{\mathbb P}

\newcommand{\QQ}{\mathbb Q}

\newcommand{\RR}{\mathbb R}

\newcommand{\X}{\mathfrak X}
\newcommand{\ZZ}{\mathbb Z}

\title{Tropical complexes}
\author{Dustin Cartwright}
\address{Department of Mathematics \\ University of Tennessee \\ 227 Ayres Hall
         Knoxville, TN 37996}
\email{cartwright@utk.edu}

\begin{document}
\begin{abstract}
We introduce tropical complexes, as an enrichment of the dual complex of a degeneration with additional data from non-transverse intersection numbers. We define cycles, divisors, and linear equivalence on tropical complexes, analogous both to the corresponding theories on algebraic varieties and to previous work on graphs and abstract tropical curves. In addition, we establish conditions for the divisor-curve intersection numbers on a tropical complex to agree with the generic fiber of a degeneration.
\end{abstract}
\maketitle

\section{Introduction}

In~\cite{baker-norine}, Baker and Norine developed a theory of divisors and
linear equivalence on finite graphs in analogy with the theory of algebraic
curves and the specialization inequality in~\cite{baker} established a formal
connection between linear equivalence on graphs and on curves. This paper
develops a higher-dimensional generalization of this theory on what we call
tropical complexes. A tropical complex consists of a $\Delta$-complex, such as
the dual complex of a degeneration, enriched with additional integers, which are
intersection numbers from the degeneration in the case of the tropical complex
of a degeneration. In a 1-dimensional degeneration, these intersection numbers
are determined by the dual complex, and so a 1-dimensional tropical complex is
just a finite graph.

The additional data on a tropical complex~$\Delta$ allows us to define a sheaf
of affine linear functions on $\Delta$, and from there we borrow formalisms from
tropical geometry and algebraic geometry in order to develop further geometric
structures on $\Delta$. Borrowing ideas from the balancing condition in tropical
geometry, we define cycles on $\Delta$, with curves and Weil divisors being the
special cases of cycles of dimension and codimension~1. In addition, a piecewise
linear function~$\phi$ on $\Delta$ is analogous to a rational function on a
variety and to it we associate a formal sum of polyhedra whose union is the set
where $\phi$ is not affine linear. We then define principal divisors to be
equivalence classes of formal sums associated to piecewise linear functions, and
Cartier divisors as equivalence classes which are locally principal. Two Cartier
divisors are linearly equivalent if their difference is principal. Finally, we
can intersect a Cartier divisor $D$ with a curve $C$ by restricting the
piecewise linear functions defining $D$ and computing the associated divisor on
$C$, which gives a well-defined formal sum of points, even when $C$ is contained
in~$D$.

The main theorem in this paper relates the intersection product between Cartier
divisors and curves outlined above
to the intersection product on an algebraic variety. To do so, we define
specialization maps, denoted~$\rho$, from divisors and curves on the general
fiber of a regular strictly semistable degeneration to divisors and curves on
the tropical complex of that degeneration. In order for the intersection theory
on the general fiber to be detectable in the combinatorics of the tropical
complex, we need to assume that the degeneration is what we called
\defi{numerically faithful}, which means that all numerical classes of curves
and divisors in each irreducible component of the special fiber are represented by
linear combinations of strata (see Definition~\ref{d:numerically-faithful} for
details).
\begin{thm}\label{t:intersect-intro}
Suppose that $\X$ is a numerically faithful, regular, strictly semistable
degeneration and $\Delta$ is its tropical complex. If $D$ is a divisor on the
generic fiber of~$\X$ and $C$ is a curve, then we have an equality of
intersection numbers:
$\deg \rho(D) \cdot \rho(C) = \deg D \cdot C$.
\end{thm}

Tropical and piecewise linear methods for computing the intersection theory
of a toric variety have a long history~\cites{fulton-sturmfels,payne,katz}.
Theorem~\ref{t:intersect-intro} partially extends this work from toric varieties
and their fans to varieties with a strictly semistable degeneration and their
dual complexes. These tropical intersection theories in $\RR^n$ depend on its
integral affine structure, and, roughly speaking, the additional data of a
tropical complex shows how the integral affine structure on one simplex extends
to adjacent simplices.

In addition to proving Theorem~\ref{t:intersect-intro}, a major purpose of this
paper is to lay the foundations for the applications of tropical complexes in
two follow-up papers. The first, \cite{cartwright-surfaces}, looks at the
combinatorial properties of $2$-dimensional tropical complexes, proving
analogues of the Hodge index theorem and Noether's formula. The second,
\cite{cartwright-specialization}, generalizes Baker's specialization inequality
for degenerations of curves~\cite{baker} to degenerations of higher-dimensional
varieties, using the definition of linear equivalence of divisors on tropical
complexes given in this paper. Some of the developments in this paper are
included to provide foundations for those results. For example, we allow
non-regular $\Delta$-complexes, which do not appear in our framework of
degenerations,
but are a natural setting combinatorially. Second, we study a
property of degenerations called robustness, which only plays a tangential role
in the main results, but is a basic concept in the specialization inequality.

We also define linear and algebraic equivalence between divisors, which is not
necessary for Theorem~\ref{t:intersect-intro}, but linear equivalence is used to
define linear series for the specialization inequality
in~\cite{cartwright-specialization} and algebraic equivalence is used in an
essential way in the proofs in~\cite{cartwright-surfaces}. Furthermore, we can
define the Picard group $\Pic(\Delta)$ and the N\'eron-Severi group
$\NS(\Delta)$ of a tropical complex~$\Delta$ to be the groups of Cartier
divisors modulo linear equivalence and algebraic equivalence respectively. Then,
we have the following structure for $\Pic(\Delta)$:
\begin{thm}\label{t:picard-structure}
The N\'eron-Severi group
$\NS(\Delta)$ of a tropical complex $\Delta$ is finitely generated. If
$\Jac(\Delta)$ denotes the kernel of the homomorphism from $\Pic(\Delta)$ to
$\NS(\Delta)$, then $\Jac(\Delta)$ is isomorphic to the quotient of
the additive group of a real vector space by a discrete subgroup.
\end{thm}
\noindent Theorem~\ref{t:picard-structure} generalizes the well-known structure
of the Picard group of a tropical curve as the product of a tropical Jacobian by
$\ZZ$, with the second factor representing the degree of the divisor.

While we have used the term ``tropical complex'' throughout the above
discussion, many of the constructions and basic results in this paper work
for a more general class, called weak tropical complexes, and when
possible, we state the results in those terms. Nonetheless,
Theorem~\ref{t:intersect-intro}, as well as the main results
of~\cites{cartwright-surfaces,cartwright-specialization} work only for tropical
complexes.

Duval, Klivans, and Martin have also introduced a generalization of linear
equivalence from graphs to higher-dimensional
complexes~\cite{duval-klivans-martin}, which differs significantly from ours.
They define a critical group for an arbitrary simplicial complex, without the
additional data which is essential in our constructions. Moreover, they define a
group in each dimension, which consists of formal sums of simplices modulo
chip-firing operations, which are indexed by simplices of the same dimension. In
contrast, for divisors supported on the codimension~$1$ skeleton of a tropical
complex, the chip-firing moves correspond to functions which are linear on each
simplex, and thus are generated by the vertices. It is only for $1$-dimensional
tropical complexes that the divisors supported on the codimension~$1$ cells and
the chip-firing moves are indexed by the same set, and the two constructions are
related.

Tropical complexes are not the only way of packaging information about a
degeneration of algebraic varieties. As in the case of
curves~\cites{amini-baker,katz-zureick-brown}, a semistable degeneration gives
rise to a semisimplicial object in the category of smooth schemes and it would
be interesting to develop a theory and applications of divisors and
intersections in such a setting. In a different direction, we work only over
discrete valuations, for which our schemes are Noetherian and there is a theory
of regular models. Accordingly, tropical complexes are intrinsically discrete
objects. Over a non-discrete, rank~$1$ valuation ring, continuous lengths would
be needed along the edges of the complex, and continuous variations would also
be necessary for any non-trivial sort of moduli space of tropical complexes,
comparable to~\cite{brannetti-melo-viviani}. 

The remainder of the paper is organized as follows. In
Section~\ref{s:complexes}, we define tropical complexes and weak tropical
complexes and explain their construction from strict semistable degenerations.
Section~\ref{s:linear} defines a sheaf of linear functions on tropical
complexes, which are used to define cycles on a weak tropical complex.
Section~\ref{s:cartier} introduces Cartier divisors and contains the proof of
Theorem~\ref{t:intersect-intro}. Section~\ref{s:equiv} introduces linear and
algebraic equivalence on divisors and contains the proof of
Theorem~\ref{t:picard-structure}.

\subsection*{Acknowledgments}

Throughout this project, I've benefited from conversations with Matt Baker,
Spencer Backman, Alex Fink, Christian Haase, Paul Hacking, June Huh, Eric Katz,
Madhusudan Manjunath, Farbod Shokrieh, Bernd Sturmfels, Yu-jong Tzeng, and
Josephine Yu. I'd especially like to thank Sam Payne for his many insightful
suggestions and thoughtful comments on an early draft of this paper. I was
supported by the National Science Foundation award number DMS-1103856 and
National Security Agency award H98230-16-1-0019. I also thank
the University of Georgia for its hospitality during the first round
of revisions on this paper.

\section{Degenerations and their tropical complexes}\label{s:complexes}

This section defines tropical complexes as combinatorial
objects and explains their construction from a degeneration. 
We begin by discussing the combinatorial
properties of degenerations.

We use \defi{degeneration} to refer to a regular scheme $\X$ which is flat and
proper over a discrete valuation ring $R$, and such that the
special fiber $\X_0$ is a reduced simple normal crossing divisor. In addition,
we always assume that the residue field of the discrete
valuation ring~$R$ is algebraically closed.

We will write $\X_\eta$ for the generic fiber of a degeneration~$\X$, and
$n$ will denote the dimension of $\X_\eta$, which we will also call the
relative dimension of $\X$. If $R$ contains a field of
characteristic~$0$, then starting from any smooth proper variety~$\X_\eta$ over
the fraction field of~$R$, it is possible to find a degeneration~$\X$, possibly
after a ramified extension, by the semistable reduction theorem of Knudsen,
Mumford, and Waterman~\cite{kkmsd}*{p.~53}.

The dual complex of the degeneration~$\X$ is a $\Delta$-complex which encodes
the combinatorics of the intersections between components of the special
fiber~$\X_0$. For details on $\Delta$-complexes, we refer
to~\cite{hatcher}*{Sec.~2.1} or to~\cite{kozlov}*{Def.~2.44}, where it is called
the ``gluing data for a triangulated space.'' In brief, the data of
a~$\Delta$-complex consists of a finite set for each $k \geq 0$, called the set
of $k$-dimensional simplices, and for each index $0 \leq i \leq k$ a face map
$d_{k,i}$ from the $k$-dimensional simplices to the $(k-1)$-dimensional
simplices, satisfying the compatibility condition that $d_{k-1,i} \circ d_{k,j}
= d_{k-1,j-1} \circ d_{k,i}$ for all $0 \leq i < j \leq k$. We will always
indicate the set of $k$-dimensional simplices with a subscript, such as
$\Delta_k$. The face map indicates how to glue the faces of a topological
simplices, from which a $\Delta$-complex specifies a topological space, its
geometric realization. In this paper, we will not usually distinguish between
the combinatorial data of a $\Delta$-complex and its geometric realization.

The dual complex~$\Delta$ of a degeneration~$\X$ has vertices (i.e.\
0-dimensional simplices) which are in bijection with the irreducible components
of~$\X_0$, and we pick an arbitrary ordering on the vertices. For a vertex~$v
\in \Delta_0$, we write $C_v$ for the corresponding component. For each
subset~$I = \{v_0, \ldots, v_k\}$ of the vertices, with $v_0 < \ldots < v_k$,
the simple normal crossing hypothesis implies that $\cap_{v \in I} C_{v}$ is a
disjoint union of smooth varieties of dimension $n - k$. We then have one
simplex~$s$ of dimension $k$ in~$\Delta$ for each of these component varieties,
and we denote the variety corresponding to $s$ by $C_s$, and call it an
\defi{$(n-k)$-dimensional stratum} of~$\X$. If we remove $v_i$ from $I$, then
$\cap_{v \in I \setminus \{v_i\}} C_v$ is a disjoint union of
$(n-k+1)$-dimensional smooth varieties and we set the face $d_{k,i}(s)$ of~$s$
to be the $(k-1)$-dimensional simplex~$s'$ such that $C_{s'}$ is the component
of this disjoint union which contains $C_s$.

For a weak tropical complex, we wish to additionally record the degree of the
intersection of each $n$-dimensional stratum $C_v$, considered as a divisor
in~$\X$ with each curve $C_r$, corresponding to an $(n-1)$-dimensional
simplex~$r$. We use \defi{ridge} to refer to the
$(n-1)$-dimensional simplices of~$\Delta$. If $v$ is
not contained in~$r$, then this
intersection is transverse and the intersection number $C_v \cdot C_r$ is equal
to the cardinality of the intersection $C_v \cap C_r$, which is equal
to the number of $n$-dimensional simplices containing both $v$ and~$r$. On the
other hand, if $v$ is contained in~$r$, then $C_v \supset C_r$ and the degree of
the intersection $C_v \cdot
C_r$ is not, in general, determined by the dual complex. We introduce
notation to record the intersection number in terms of $v$ and~$r$:
\begin{equation*}
\alpha(v,r) = - \deg C_v \cdot C_r,
\end{equation*}
We refer to the integers $\alpha(v, r)$ as \defi{structure constants}, and the
combination of $\Delta$ and the structure constants as the \defi{weak tropical
complex of $\X$}.

When the relative dimension of~$\X$ is~$1$, then a component $C_v$ is a
curve, and so the only case when the intersection of $C_v$ with $C_r$ is not
transverse is when $r = v$. In this case, the self-intersection number of~$C_v$
is determined by the dual complex, and $\alpha(v, v)$ is the degree of $v$
in~$\Delta$~\cite{liu}*{Prop.~9.1.21b}. This generalizes to a constraint on the
structure constants in higher dimensions, but it does not determine them fully:

\begin{prop}\label{p:constraint}
If $\Delta$ and $\alpha$ are the weak tropical complex of a degeneration~$\X$,
and $r$ is any ridge of $\Delta$,
then:
\begin{equation}\label{eq:constraint-in-motivation}
\sum_{v \in r_0} \alpha(v, r) = \deg r,
\end{equation}
where the summation is over the vertices~$v$ of $r$, and $\deg r$ is number of
$n$-dimensional simplices containing $r$.
\end{prop}

\begin{proof}
Since the special fiber $\X_0$ is reduced, the divisor associated to a
uniformizer $\pi$ of $R$ is the sum of irreducible divisors $\sum_{v \in
\Delta_0} [C_v]$ as $v$ ranges over all vertices of~$\Delta$. We consider the
intersection of the special fiber with the curve $C_r$ corresponding
to the fixed ridge~$r$. If $v$ is a vertex
of~$r$, then the intersection of $C_v$ has degree $-\alpha(v,r)$ by the
definition of the structure constants. On the other hand, if $v$ is not
contained in $r$, then the intersection of $C_v$ with $C_r$ is the disjoint
union of reduced points whose total cardinality is the number of simplices
containing both $v$ and $r$. By linearity, the intersection of
$\X_0 = \sum_{v \in \Delta} [C_v]$ with $C_r$ has degree:
\begin{equation*}
\sum_{v \in r_0} -\alpha(v, r) + \deg r = 0,
\end{equation*}
because the $\X_0$ is a principal divisor.
Rearranging, we've verified~(\ref{eq:constraint-in-motivation}).
\end{proof}

When $n$ is at least $2$, we have the following interpretation of the structure
constants:

\begin{figure}
\includegraphics{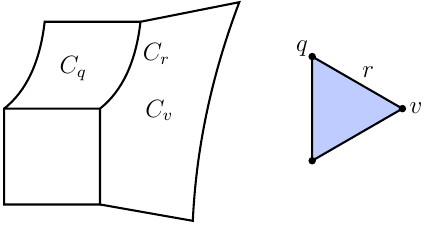}
\caption{On the left is the special fiber of a degeneration~$\X$ of relative
dimension $n=2$, consisting of three components. On the right is the dual
complex of this degeneration. In this case,
Lemma~\ref{l:alpha-self-intersection} states that $\alpha(v,r)$, which is
defined as $-\deg C_v \cdot C_r$, where the intersection is taken in~$\X$, is
equal to $-\deg C_r^2$, taken in the surface
$C_q$.}\label{f:self_intersection}
\end{figure}
\begin{lem}\label{l:alpha-self-intersection}
Suppose that $\X$ is a degeneration of relative dimension at
least~$2$, with dual complex $\Delta$. Let $v$ and $r$ be a vertex and ridge,
respectively, of~$\Delta$. Let $q$ be the unique
$(n-2)$-dimensional face of $r$ which does not contain $v$, as illustrated in
Figure~\ref{f:self_intersection}. Then:
\begin{equation*}
\alpha(v,r) = -\deg C_r^2,
\end{equation*}
where the self-intersection of the curve $C_r$ is taken in the surface~$C_q$.
\end{lem}

\begin{proof}
By definition, $-\alpha(v,r)$ is the degree of the restriction of the Cartier
divisor $C_v$ to the curve~$C_r$. Let $q$ be as in the statement. We then
restrict the Cartier divisor $C_v$ first to $C_q$ and then to $C_r$. Since $q$
does not contain $v$, then the intersection of $C_v$ with $C_q$ is transverse
and consists of a disjoint union of smooth curves, including $C_r$ as well as
$C_{r'}$ for any other ridge $r'$ containing both $v$ and $q$. Since these
curves are disjoint, those other than $C_r$ do not affect the pullback to $C_r$,
and so and so the pullback of $C_v$ from $\X$ to $C_r$ is the same as the
pullback of the divisor $C_r$ from $C_q$ to $C_r$, which is the
self-intersection of $C_r$.
\end{proof}

Now suppose we fix an $(n-2)$-dimensional simplex~$q$, and look at the
corresponding surface $C_q$, which contains the $1$-dimensional stratum~$C_r$
for each ridge~$r$ containing~$q$. Any two distinct curves $C_r$ and $C_{r'}$
intersect at a finite union of reduced points, which are in bijection with
simplices containing both $r$ and $r'$. Likewise, by
Lemma~\ref{l:alpha-self-intersection}, the self-intersection of a single curve
$C_r$ is $-\alpha(\opp_q(r), r)$, where $\opp_q(r)$ is defined to be the vertex
of $r$ not contained in~$q$. Therefore, the intersection pairing of
the curves $C_r$ on the surface $C_q$ is determined by the topology
of $\Delta$ and the structure constants $\alpha(v,r)$. More precisely, the
matrix $M_q$
representing this pairing is the matrix
whose rows and columns are labeled by the ridges $r$ containing~$q$ and whose
entries are:
\begin{equation}\label{eq:local-intersection-matrix-realizable}
(M_q)_{r, r'} = C_r \cdot C_{r'} =
\begin{cases}
\# \{\mbox{simplices containing $r$ and $r'$}\} & \mbox{if } r \neq r' \\
-\alpha(\opp_q(r), r) & \mbox{if } r = r'
\end{cases}
\end{equation}
As a consequence of the Hodge index theorem, we have:

\begin{prop}\label{p:at-most-one-positive}
Given a degeneration $\X$ and an $(n-2)$-dimensional simplex~$q$ of the dual
complex $\Delta$, the matrix $M_q$ defined in
(\ref{eq:local-intersection-matrix-realizable}) has at most 1 positive
eigenvalue.
\end{prop}

\begin{proof}
By Lemma~\ref{l:alpha-self-intersection}, the matrix $M_q$ contains the
intersection numbers between the $1$-dimensional strata in the surface $C_q$.
Since $C_q$ is a smooth, proper surface, it is projective. Therefore,
the Hodge index theorem implies that the intersection matrix between divisors
which form a generating set for all numerical classes on $C_q$ has exactly one
positive eigenvalue. Since $M_q$ restricts this pairing to a subspace, it has at
most one positive eigenvalue. 
\end{proof}

\begin{ex}\label{ex:k3-non-robust}
Consider the family $\X' \subset \PP^3_{\CC[[t]]}$ defined by the
equation $xyzw + t(x^4 + y^4 + z^4 + w^4)$, where $x$, $y$, $z$, and~$w$ are the
coordinates of~$\PP^3$. The special fiber consists of $4$ planes meeting along
their coordinate lines. However, $\X'$ is not regular, with 24 singularities at the points with coordinates $(0:0:1:\zeta)$ where $\zeta$ is a
primitive $8$th root of unity, together with permutations of these coordinates.

These singularities are ordinary double point singularities. For example, we can
look at the local equation near the point $(0:0:1:\zeta)$ by substituting
$z = 1$ and $w = w' + \zeta$ into the defining equation to get:
\begin{equation*}
\zeta xy + 4 \zeta^3 t w' + \mbox{higher order terms}.
\end{equation*}
The degree 2 part of this equation defines a non-degenerate quadratic form,
which is a characterization of an ordinary double point singularity. An ordinary
double point singularity on a 3-fold has two different small resolutions, in the
sense that they don't add any new divisors.
In our context, these two resolutions 
can be obtained by blowing up one of two coordinate planes in the special
fiber passing through the singular point. The strict transform of the chosen
coordinate plane is isomorphic to the blow-up of a plane at all the singular
points, but the strict transforms of the other planes are unaffected by the
blow-up.

We let $\X$ be the family formed by successively blowing up $\widetilde\X$ at
the components of the special fiber, yielding a degeneration. Then, as
with $\widetilde\X$, the special fiber of $\X$ consists of 4 components, which
therefore correspond to four vertices in the dual complex, and we denote these
vertices $v_1$, $v_2$, $v_3$, and $v_4$, in the same order as the blow-ups.
Therefore, the surface $C_{v_1}$ consists of the plane $\PP^2$ blown up at 4
points along each of the 3 coordinate lines, $C_{v_2}$ is $\PP^2$ blown up at 4
points along each of the 2 coordinate lines which don't intersect
$C_{v_1}$, and so on. Each pair of these components meets along the strict
transform of a coordinate line, and any three of them meet at the strict
transform of a coordinate point in $\PP^3_{\CC}$. Therefore, the dual complex of
$\X$ is the boundary of a tetrahedron. We denote the edges $e_{12}$, $e_{13}$,
and so on.

Since blowing up a point decreases the self-intersection of a curve in a surface
by $1$, the self-intersection of a coordinate line $C_{e_{ij}}$ in one of a component
$C_{v_i}$ is
either $1$, if $i < j$ and so $C_{v_i}$ was not blown up at points in
$C_{e_{ij}}$, or $-3$ if $i > j$ and so it was.
Therefore, by Lemma~\ref{l:alpha-self-intersection}, 
\begin{equation*}
\alpha(v_j, e_{ij}) =
\begin{cases}
-1 & \mbox{if } i < j \\
3 & \mbox{if } i > j.
\end{cases}
\end{equation*}
One immediately verifies Proposition~\ref{p:constraint} since the degree of each
edge is $2$.

From the structure constants, we can compute the intersection matrices as
in~(\ref{eq:local-intersection-matrix-realizable}) in order to verify
Proposition~\ref{p:at-most-one-positive}:
\begin{gather*}
M_{v_1} = \begin{bmatrix}
-3 & 1 & 1 \\
1 & -3 & 1 \\
1 & 1 & -3 \\
\end{bmatrix}
M_{v_2} = \begin{bmatrix}
1 & 1 & 1 \\
1 & -3 & 1 \\
1 & 1 & -3 \\
\end{bmatrix} \\
M_{v_3} = \begin{bmatrix}
1 & 1 & 1 \\
1 & 1 & 1 \\
1 & 1 & -3 \\
\end{bmatrix}
M_{v_4} = \begin{bmatrix}
1 & 1 & 1 \\
1 & 1 & 1 \\
1 & 1 & 1 \\
\end{bmatrix}
\end{gather*}
Of these, $M_{v_1}$ is negative definite, and so has no positive eigenvalues,
while the other three matrices each have exactly one positive eigenvalue.

See also Example 1.3 in \cite{yu-balancing} for another discussion of the
combinatorics of this same degeneration.
\end{ex}

While Example~\ref{ex:k3-non-robust} shows that the intersection matrix
can have either $0$ or~$1$ positive eigenvalues, we will be most interested in
degenerations where every intersection matrix falls in the latter case. We
want the weak tropical complex to reflect the algebraic geometry of the
degeneration, and capturing the Hodge index theorem for surfaces is a key part
of that geometry, as shown in Proposition~\ref{p:specialization-divisor} and in
the follow-up papers~\cites{cartwright-surfaces,cartwright-specialization}.

We wish to define tropical complexes as abstract objects, consisting of a
$\Delta$-complex and structure constants $\alpha(v,r)$, not necessarily coming
from a degeneration, but following the conclusions of
Propositions~\ref{p:constraint} and~\ref{p:at-most-one-positive}. One way in
which a tropical complex will be more general than the data coming from a
degeneration is that we allow non-regular $\Delta$-complexes. A regular
$\Delta$-complex is one in which the faces of a single simplex are always
distinct~\cite{kozlov}*{Def.~2.47}, whereas a non-regular $\Delta$-complex
allows distinct faces of a single simplex to be identified with each other. As a
combinatorial structure, there is no need to restrict to regular
$\Delta$-complexes, and this paper provides foundations for the further
combinatorial theory of surfaces in~\cite{cartwright-surfaces}, which does not
require regularity. In this paper, we will use a non-regular $\Delta$-complex to
illustrate algebraic equivalence and the structure of the Picard group in
Example~\ref{ex:torus}. We could subdivide each simplex to obtain a regular
$\Delta$-complex, but the advantage of the non-regular $\Delta$-complex is that
it requires only $2$ maximal simplices, rather than $8$. Moreover, while the
dual complex of a simple normal crossing divisor is always a regular
$\Delta$-complex, non-regular $\Delta$-complexes naturally arise from the
regular semistable degenerations with self-intersections, at least as long as
the self-intersections have no monodromy, as defined
in~\cite{abramovich-caporaso-payne}*{Def.~6.2.2}.

Since the faces of a simplex in a non-regular $\Delta$-complex may have been
identified with each other, a given simplex~$s$ in the geometric realization of
a $\Delta$-complex~$S$ may not be homeomorphic to a standard simplex. Thus, we
will use the term
\defi{parameterizing simplex}~$\tilde s$, to denote the standard simplex, before
such identifications, which maps surjectively onto $s$. In particular, the
parameterizing simplex of an $k$-dimensional simplex always has exactly $k+1$
vertices, some of which may be identified in $S$.

The local topology of a $\Delta$-complex $S$ around any point in the interior of
a $k$-dimensional simplex $s$ can be represented by a $\Delta$-complex called
its \defi{link} and denoted $\link_S(s)$. In the case of a regular
$\Delta$-complex
$S$, the $m$-dimensional simplices of $\link_S(s)$ are in bijection with the
$(m+k+1)$-dimensional simplices of~$S$ which contain $s$, and this
bijection preserves the containment relation. In the case of a non-regular
$\Delta$-complex~$S$, the $m$-dimensional simplices of $\link_S(s)$ are in
bijection with pairs of an $(m+k+1)$-dimenisonal simplex~$s'$ of $S$ together
with an identification of $s$ as a face of $s'$. For example, if $S$ is a graph
with a loop at a single vertex $v$, then $\link_v(S)$ consists of two vertices,
corresponding to the two different ways that $v$ is a
face of the unique edge of~$S$. The geometric realization of $\link_v(S)$ is two
disjoint points, which is the intersection of $S$ with a small ball around $v$
in a (reasonable) embedding of $S$ in Euclidean space. For a precise definition and further details on the link of a
$\Delta$-complex, see~\cite{kozlov}*{p.~31}.

Now suppose that $s$ is a $k$-dimensional simplex of~$S$ and $t$ is a
vertex in $\link_S(s)$. Since $t$ corresponds to an identification of $s$
with a face of a $(k+1)$-dimensional simplex $s'$, there exists a
unique vertex of the parameterizing simplex $\tilde s'$ not contained in the
face identified with $t$, and we denote this vertex by $\opp_s(t)$.

\begin{defn}\label{d:weak-tropical-complex}
An \defi{$n$-dimensional weak tropical complex} is a pair $(S, \alpha)$, where
$S$ is a finite, connected $\Delta$-complex of dimension at most $n$, and $\alpha$ is an
integer-valued function on the set of pairs $(v,r)$ such that $r$ is an
$(n-1)$-dimensional simplex of $S$ and $v$ is a vertex of the parameterizing
simplex~$\tilde r$,
such that, for each ridge $r$, \begin{equation}\label{eq:constraint}
\sum_{v \in \tilde r} \alpha(v, r) = \deg(r),
\end{equation}
where the summation is over the vertices of the parameterizing simplex of~$r$.
Here, $\deg(r)$ refers to the cardinality of the finite set $\link_S(r)$. 
The values $\alpha(v,r)$ are called the \defi{structure
constants} of~$\Delta$.

We will refer to the $n$-dimensional and $(n-1)$-dimensional simplices of~$S$ as
the \defi{facets} and \defi{ridges} of $\Delta$, respectively. Moreover, we
will, in general, not distinguish notationally between the weak tropical complex
$\Delta$ and the underlying $\Delta$-complex $S$. Thus, we will also refer to
vertices or points
to refer to vertices or points of (the geometric realization of) $S$, and
$\Delta_k$ will refer to the set of $k$-dimensional simplices of $S$.
\end{defn}

Definition~\ref{d:weak-tropical-complex} is consistent with our previous
terminology of the weak tropical complex of a degeneration:

\begin{prop}\label{p:weak-tropical-complex-degeneration}
The weak tropical complex of a degeneration of relative dimension $n$ satisfies
the conditions in Definition~\ref{d:weak-tropical-complex}.
\end{prop}

\begin{proof}
This follows immediately from Proposition~\ref{p:constraint}.
\end{proof}

\begin{rmk}
While the simplices of an $n$-dimensional simplicial complex all have dimension
at most~$n$, this bound does not have to be achieved. The maximal simplices of
the underlying $\Delta$-complex may have different dimensions, or all have
dimension less than~$n$. Thus, the dimension of a weak tropical complex is part
of the data of the weak tropical complex and not a property of the underlying
$\Delta$-complex. In this paper, $n$ will always denote the dimension of a weak
tropical complex.
\end{rmk}

\begin{defn}\label{d:tropical-complex}
Let $\Delta = (S, \alpha)$ be an $n$-dimensional weak tropical complex, and fix
an $(n-2)$-dimensional simplex~$q$ of~$\Delta$. We define the \defi{local
intersection matrix} at $q$
to be the symmetric matrix $M_q$ whose rows and columns are labeled by the
vertices of $\link_\Delta(q)$, and such that the entry whose row and column
correspond to $t$ and $u$, respectively, in $\link_\Delta(q)$ is:
\begin{equation}\label{eq:local-intersection-matrix}
(M_q)_{t,u} = \begin{cases}
\# \{\mbox{edges between $t$ and $u$ in $\link_\Delta(q)$}\} &
  \mbox{if $t \neq u$} \\
-\alpha(\opp_q(t),r(t)) + 2 \cdot\# \{\mbox{loops at $t$ in $\link_\Delta(q)$}\}&
  \mbox{if $t = u$}
\end{cases}
\end{equation}

An \defi{$n$-dimensional tropical complex} is an $n$-dimensional weak tropical
complex~$\Delta$ such that, for each $(n-2)$-dimensional simplex $q$,
of~$\Delta$, the local intersection matrix $M_q$ has exactly one positive
eigenvalue.
\end{defn}

If $\Delta$ is the weak tropical complex of a degeneration, then it is regular,
and so $\link_\Delta(q)$ contains no loops, and the local intersection
matrix~$M_q$ defined in~(\ref{eq:local-intersection-matrix}) agrees with the
intersection matrix of the surface~$C_q$ given
in~(\ref{eq:local-intersection-matrix-realizable}). The latter matrix is the
adjacency matrix of $\link_\Delta(q)$, with structure constants along the
diagonal. When $\link_\Delta(q)$ has loops, one convention for its
adjacency matrix is that a diagonal entry is equal to twice the number of loops
at the corresponding vertex,
which preserves the property that the sum of a row is the degree of
the corresponding vertex.
With this convention, $M_q$ is the difference
of the adjacency matrix of the graph $\link_\Delta(q)$ and a diagonal matrix
containing structure constants, which motivates the second term of the
diagonal case of~(\ref{eq:local-intersection-matrix}).

Example~\ref{ex:k3-non-robust} shows that the weak tropical complex of a 
degeneration is not necessarily a tropical complex. To obtain a tropical
complex, we need to put a condition on the degeneration. Recall that a Cartier
divisor on a normal variety is called \defi{big} if, for some multiple of the
divisor, the rational map to projective space defined by taking the complete
linear series is birational onto its image~\cite{lazarsfeld}*{Sec. 2.2}.

\begin{defn}\label{d:robust}
Let $\X$ be a degeneration of relative dimension~$n$ with $\Delta$ its
dual complex. If $s$ is a simplex of $\Delta$ of dimension $n-k$, we write $D_s$
for the divisor on $C_s$ given as the sum $\sum [C_{s'}]$  over all
$(n-k+1)$-dimensional simplices $s'$ containing $s$.

We say that $\X$ is \defi{robust at $s$} if $D_s$ is a big divisor on $C_s$. If
$k$ is an integer in the range $1 \leq k \leq n$, we say that $\X$ is robust in
dimension $k$ if for each simplex $s$ of dimension $n-k$, $\X$ is robust at~$s$.
\end{defn}

\begin{prop}\label{p:tropical-complex-robustness}
Suppose that $\Delta$ is the weak tropical complex of a 
degeneration~$\X$. Then, $\Delta$ is a tropical complex if and only if $\X$ is
robust in dimension~$2$.
\end{prop}

\begin{proof}
We first suppose that $\X$ is robust in dimension~$2$ and fix an
$(n-2)$-dimensional simplex $q$ of~$\Delta$. As in the discussion before
Proposition~\ref{p:at-most-one-positive}, $M_q$ is the intersection matrix
on~$C_q$ of the components of~$D_q$, and by
Proposition~\ref{p:at-most-one-positive}, it has at most one positive
eigenvalue. Therefore, it will suffice to find a linear combination of the
components which has positive self-intersection. By assumption, $D_q \subset
C_q$ is a big divisor. We take a sufficiently large multiple of~$D_q$ such that
it defines a rational map to $\PP^N$ with two-dimensional image. If we remove
any divisors in the base locus of this map, we will get a divisor~$D$ such that
$D^2 > 0$. Since the base locus is contained in~$D_q$, then $D$ is a linear
combination of the components of~$D_q$. Thus, together with
Proposition~\ref{p:at-most-one-positive}, we've shown that $M_q$ has one
positive eigenvalue.

Conversely, suppose that $\Delta$ is a tropical complex. Then, by definition,
for each $(n-2)$-dimensional simplex~$q$, the local intersection matrix~$M_q$
has at least one positive eigenvalue. By rationally approximating the
corresponding eigenvector and then scaling, we can then find an integer
vector~$v$ such that $v^T M_q v$ is positive. Using the entries of this vector
as the coefficients for a linear combination of the components of~$D_q$, we get
a divisor~$D$ on $C_q$ with positive self-intersection. Let $A$ be an ample
divisor on~$C_q$. By the Hodge index theorem, $A \cdot D$ is non-zero, and if
necessary, we replace $D$ with $-D$, so that $A \cdot D$ has positive degree.
Then, for $m$ sufficiently large, $K_{C_q}-mD$ will have negative degree with
respect to~$A$ and thus not be linearly equivalent to any effective divisor.
Therefore, Riemann-Roch formula for surfaces implies that $h^0(C_q, mD)$ grows
quadratically in~$m$, so $D$ is big by~\cite{lazarsfeld}*{Lem.~2.2.3}. If $k$ is
the largest coefficient of a component of~$D$, then $kD_q$ is the sum of $D$
with an effective divisor, so $kD_q$ is big, and thus $D_q$. Therefore, $\X$ is
robust in dimension~$2$.
\end{proof}

\begin{ex}\label{ex:k3-robust}
As in Example~\ref{ex:k3-non-robust}, we start with the family $\X'
\subset \PP^3_{\CC[[t]]}$ defined by $xyzw + t(x^4 + y^4 + z^4 + w^4)$, and each
singularity of $\X'$ is an ordinary double point singularity.
An ordinary double point singularity on a 3-fold admits two different small resolutions, both of
which introduce a rational curve over the singular point, which is contained in
only one of the two components of $\X_0$ containing the point. Unlike
Example~\ref{ex:k3-non-robust}, where we resolved by blowing up the components
of $\X_0$, here we resolve each of the 24 singularities of $\X'$ individually,
such that on each intersection of two components, the four singularities are
resolved symmetrically, in the sense that two of the exceptional curves are
contained in one component and two are contained in the other.

Therefore, if $\X$ denotes the resulting degeneration, then the special fiber
$\X_0$ consists of four components, each isomorphic to $\PP^2$ blown up at six
points. The divisor $D_v$ on each of these components consists of the strict
transform of the union of three lines, each passing through two of the points
that were blown up. Thus, each line has self-intersection $-1$, so the weak
tropical complex of $\X$ has $\alpha(v,e) = 1$ for every vertex $v$ of every
edge $e$. One can check that these structure constants define a tropical
complex, and also that $D_v$ is a big divisor on $C_v$, making $\X$ a robust
degeneration.
\end{ex}

Robustness in dimensions greater than 2 will play a critical role in the
specialization theorem for tropical complexes~\cite{cartwright-specialization},
but will not be important for the main results of this paper. Instead, the
hypothesis that Theorem~\ref{t:intersect-intro} puts on a degeneration is that
it is numerically faithful, which we now explain. Recall that two divisors
(resp.\ curves) on a smooth variety are \defi{numerically equivalent} if they
have their intersection numbers are equal, when paired with any curve (resp.\
divisor) on the variety~\cite{lazarsfeld}*{Defs. 1.1.14 and 1.4.25}.

\begin{defn}\label{d:numerically-faithful}
We say that a degeneration~$\X$ is \defi{numerically faithful at a
vertex~$v$} of a dual complex~$\Delta$ if for every divisor~$D$ in~$C_v$, some
multiple of~$D$ is numerically equivalent to a linear combination of the
$(n-1)$-dimensional strata contained in $C_v$, and for every curve
$C$ in $C_v$, some multiple of $C$ is numerically equivalent to a linear
combination of the $1$-dimensional strata contained in~$C_v$.

A degeneration~$\X$ is \defi{numerically faithful} if it is numerically
faithful at every vertex of $\Delta$.
\end{defn}

\begin{prop}\label{p:numerically-faithful-robust}
Let $\X$ be a degeneration which is numerically faithful at a vertex $v$ of its
weak tropical complex $\Delta$ and suppose that $C_v$ is projective. Then
$\Delta$ is robust at~$v$.
\end{prop}

\begin{proof}
Because $C_v$ has a big divisor, and big divisors on a projective variety are
characterized by their numerical equivalence
class~\cite{lazarsfeld}*{Cor.~2.2.8}, there is a linear combination of the
$(n-1)$-dimensional strata on $C_v$ which is big. We can then add an effective
divisor so that all the $(n-1)$-dimensional strata have the same coefficient,
and so a multiple of $D_v$ is big, and thus $D_v$ is big.
\end{proof}

The degeneration in Example~\ref{ex:k3-robust} shows that
the converse to Proposition~\ref{p:numerically-faithful-robust} is not true,
because it is robust, but not numerically faithful.
Each
surface~$C_v$ of the special fiber in that example is the blow-up of $\PP^2$ at
$6$ points, which has Picard group~$\ZZ^7$ with a non-degenerate intersection
product, but there are only three boundary curves~$C_e$ contained in each $C_v$,
so they do not span the numerical equivalence classes.

On the other hand, the following proposition shows that degenerations where the
special fiber is built from hyperplane complements are always numerically
faithful. Such degenerations are related to tropical manifolds, which are defined to
to have the local structure of the Bergman fan of a
matroid~\cites{mikhalkin,jell-rau-shaw,jell-shaw-smacka} and are used to compute
the homology of the general fiber in~\cite{ikmz}.

\begin{prop}\label{p:hyperplane-complement}
Let $\X$ be a degeneration. Suppose that, for each $k$-dimensional simplex
$s \in \Delta_k$, where $k \leq n-1$, the complement $C_s \setminus D_s$ is
isomorphic to $\PP^{n-k} \setminus H$, where $H$ is the union of one or more
hyperplanes. Then $\X$ is numerically faithful.
\end{prop}

\begin{proof}
Let $C_s$ be the stratum corresponding to a simplex $s$ of dimension~$k$
in~$\Delta$, and let $Z$ be a cycle of dimension $m$. By assumption,
$C_s \setminus D_s$ is isomorphic to an open subset of $\mathbb A^{n-k}$, whose
Chow groups are trivial in all dimensions less than $n-k$. Thus, if $m$
is less than $n-k$, then $Z$ is linearly equivalent to a
cycle contained in $D_s$. Therefore, by induction, $Z$ is linearly equivalent to a
linear combination of strata. Taking just the cases when $m$ is either $1$ or
$n-1$, this is sufficient to show that $\X$ is numerically faithful.
\end{proof}

\section{Linear functions and cycles}\label{s:linear}

In this section, we use the structure constants of a weak tropical complex to
define a sheaf of linear functions on the underlying topological space. We then
use linear functions to define a balancing condition which characterizes
cycles on a weak tropical complex.

\begin{defn}\label{d:pl-function}
A \defi{PL function} (or piecewise linear function) on an open subset of a weak
tropical complex~$\Delta$ is a continuous function $\phi$ such that on each
simplex of $\Delta$, the restriction function $\phi$ is a piecewise linear
function, with finitely many domains of linearity, and with integral slopes,
under the
identification of the simplex with the standard unit simplex.
\end{defn}

In Definition~\ref{d:pl-function} and later in this paper, we use \defi{standard
unit simplex (of dimension $k$)} to refer to the convex hull of the origin and
the unit coordinate vectors $e_1, \ldots, e_{k}$ in $\RR^{k}$. Although the
equation of a linear function will depend on the particular identification we
choose, the set of linear functions with integral slopes does not.

In addition, we define linear functions on~$\Delta$.
Combinatorially, the main difficulty is
understanding linearity across the boundaries between two facets. The following
gives a function from a neighborhood of each ridge in a weak tropical complex to
a vector space, and then we obtain all linear functions by composing with an
affine linear function on that vector space.

\begin{const}\label{con:linear-function}
Let $\Delta$ be an $n$-dimensional weak tropical complex and $r$ a ridge
of~$\Delta$. Let $N_r$ be the simplicial complex obtained by attaching one
$n$-dimensional simplex for each vertex in $\link_\Delta(r)$ onto a central
$(n-1)$-dimensional simplex, which we identify with $\tilde r$. Thus, there
exists a natural map $\pi_r \colon N_r \rightarrow \Delta$ extending the
parameterization of $r$ by by the central simplex. If $N_r^o$ denotes the open
subset of~$N_r$ consisting of the interior of the central $(n-1)$-dimensional
simplex and the interiors of the
facets of~$N_r$, then $\pi_r \vert_{N^o_r}$ is a local homeomorphism, and,
moreover, if $\Delta$ is a regular $\Delta$-complex, then $\pi_r \vert_{N^o_r}$
is an open immersion.

Now let $L_r$ be the quotient $\RR^{d + n} / \RR$ by the line generated by
the vector $(1, \ldots, 1, -\alpha(v_1, r), \ldots, -\alpha(v_n, r))$ where $d$
is the cardinality of $\link_\Delta(r)$ and $v_1,
\ldots, v_n$ denote the vertices of~$\tilde r$. Let $e_i \in L_r$ denote the
image of the $i$th coordinate vector of $\RR^{d+n}$, and let $w_1, \ldots, w_d$
denote the vertices of $N_r$ that are not contained in $\tilde r$. Then, $\phi_r\colon N_r
\rightarrow L_r$ is the map which is linear on simplices and sends $w_i$ to
$e_i$ for $1 \leq
i \leq d$ and sends $v_i \in \tilde r$ to $e_{d+ i}$.

A \defi{linear function} on an open subset $U$ of a weak tropical
complex~$\Delta$ is a PL function $\phi$ such that on each facet~$f$ meeting $U$,
$\phi \vert_{f \cap U}$ is linear, under the identification of $f$ with a unit
simplex in $\RR^n$, and for each ridge $r$ meeting $U$, the pullback $\phi \circ
\pi_r \vert_{\pi_r^{-1}(U)}$ coincides with $\ell \circ \phi_r
\lvert_{\pi_r^{-1}(U)}$, where
$\ell$ is an affine linear function
$\ell \colon L_r \rightarrow \RR$.
\end{const}

In Construction~\ref{con:linear-function}, the image $\phi_r(N_r)$
is contained the affine hyperplane in~$L_r$ consisting of vectors of the form
$\sum_{i=1}^{d+n} c_i e_i$ such that $\sum_{i=1}^{d+n} c_i = 1$. Therefore, although
there is one real parameter for the constant of an affine linear function
on~$L_r$  and $d+n-1$ integral parameters for the slopes, there are only $d+n-2$
parameters for the slope of a linear function on $N_r^o$.

For curves, linear functions are the same as those used with tropical curves,
such as in \cite{mikhalkin-zharkov}*{Def.~3.7}. The function $\phi_r$ in
Construction~\ref{con:linear-function} and its image are equivalent to the local
charts \cite{mikhalkin-zharkov}*{Def.~3.1} and the local curve model
in~\cite{mikhalkin-zharkov}*{Ex.~2.7}, respectively.

\begin{rmk}\label{r:yu-linear}
Yu has also given a definition of linear functions on the dual complex of a
degeneration in~\cite{yu-compactness}*{Def.~3.2}. Yu's definition involves not
just the data contained in a weak tropical complex, but also the group of all
numerical curve classes on each component of the special fiber. In contrast, our
definition only incorporates the numerical classes of the $1$-dimensional strata
of the degeneration, via the structure constants $\alpha(v,r)$. When the
degeneration is numerically faithful, these strata generate all numerical
classes, and the two definitions agree. However, for arbitrary degenerations,
linear functions in Yu's sense in are linear in the sense of
Construction~\ref{con:linear-function}, but not necessarily conversely.

In some cases, it may be tractable to compute the intersection numbers with
1-dimensional strata which define the weak tropical complex, but hard to
understand the entire group of numerical classes of curves. Moreover, the
structure of a weak tropical complex is sufficient for defining linear
equivalence of divisors, described in Section~\ref{s:linear} below, and the
specialization inequality~\cite{cartwright-specialization}, which builds on
linear equivalence.
\end{rmk}

\begin{rmk}\label{r:gross-siebert}
If a tropical complex happens to be homeomorphic to a manifold, then the sheaf
of linear functions gives the manifold an integral affine structure away from
the codimension~$2$ simplices. Manifolds with integral affine structures have
also been constructed from degenerations of Calabi-Yau varieties by Gross and
Siebert~\cite{gross-siebert}, building on ideas of Kontsevich and
Soibelman~\cite{kontsevich-soibelman}, but their constructions differ from ours.
Rather than regular semistable models, Gross and Siebert use what they call
toric degenerations of a Calabi-Yau variety, for which the total family~$\X$ is
allowed to have singularities, but components of the special fiber are required
to be toric varieties (see~\cite{gross-siebert}*{Def.~4.1} for the precise
definition).

In the case of Example~\ref{ex:k3-robust}, their toric degeneration would be the
family $\X' \subset \PP^3_R$ with $24$ singularities, before the
blow-ups were used to obtain a regular semistable model. In both the toric and
the regular semistable degeneration, the dual intersection complex is the
boundary of a tetrahedron, but for Gross-Siebert the singularities in the affine
linear structure lie at the midpoints of the $6$~edges, corresponding to the
singularities of~$X$~\cite{gross-siebert}*{p.~172}, whereas
Construction~\ref{con:linear-function} only puts singularities on the
$4$~vertices.

On the other hand, the sheaf of linear functions on a tropical complex agrees
with the one constructed in \cite{ghk}*{following Def.~1.4} for degenerations of
Calabi-Yau surfaces.
\end{rmk}

We now use the linear functions on a weak tropical complex to define cycles,
which are equivalences classes of formal sums of polyhedra satisfying a
balancing condition. By $k$-dimensional polyhedron in~$\Delta$, we mean a subset
of a single simplex~$s$, which is a convex polyhedron with rational slopes,
under the usual identification of $s$ with the standard unit simplex. We
sometimes work in an open subset
subset $U \subset \Delta$, in which case a $k$-dimensional polyhedron means the intersection
of a $k$-dimensional polyhedron of $\Delta$ with~$U$. A \defi{formal sum of
polyhedra} on $U$ is a finite integral linear combination of polyhedra. In
addition, we define two formal sums of $k$-dimensional polyhedra to be
\defi{equivalent} if they differ by an element of the subgroup generated by all formal
sums $[P] - [Q_1] - \cdots - [Q_r]$, where $P$ is any $k$-dimensional
polyhedron, and $Q_1, \ldots Q_r$ subdivide $P$, in the sense that $P
= Q_1 \cup \cdots \cup Q_r$ and $Q_i \cap Q_j$ is either empty or a proper face
of each.

In our constructions below, it will be useful to be able to replace a formal sum
of polyhedra with an equivalent one by subdividing sufficiently. For example,
any formal sum of polyhedra is equivalent to one where the terms form the
maximal cells of a polyhedral complex, meaning that the intersection of two
terms is either empty or a face of both.

Our goal is to define the multiplicity of a PL function $\phi$ along a formal
sum of polyhedra~$Z$, which, roughly, measures the extent to which $\phi$ is not
a linear function on $Z$. In the theory of abstract tropical curves, the divisor
of a PL function $\phi$ is defined in terms of multiplicities at each point,
which are computed as the sum of the slopes of~$\phi$ in the outgoing
directions~\cite{gathmann-kerber}*{Def. 1.4} and~\cite{baker}*{Sec.~1D}. The
number of outgoing directions is always finite, and in particular, equal to the
degree at a vertex and equal to 2 at a point in the interior of an edge. For
formal sums of polyhedra of dimension $k > 1$, the number of outgoing directions
is obviously infinite. In order to have a well-defined slope, we
limit ourselves to functions which are constant along a fixed
$(k-1)$-dimensional face, so that we can define an equivalence class of outgoing
directions, all of which will yield the same slope.

We use the lattice structure of $\ZZ^n \subset \RR^n$ to normalize the slope,
which also occurs in the
definition of the balancing condition in tropical
geometry~\cite{maclagan-sturmfels}*{Def~3.3.1}. The balancing
condition for a 1-dimensional polyhedral complex uses a minimal integral vector
in the direction of each edge, similar to the outgoing directions above. Then,
for $k$-dimensional complexes, the balancing condition is defined by first
taking the quotient by the linear span of each $(k-1)$-dimensional cell and then
applying the 1-dimensional dimension to the quotient, which inherits an integral
structure from the fact that the linear span is a rational vector
space. We use an analogous quotient to
define the slope of a function on a polyhedron.

\begin{figure}
\includegraphics{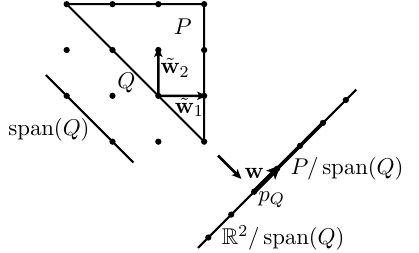}
\caption{An illustration of the quotient used in Construction~\ref{con:slope}
for defining the slope. Here, $P$ is a $2$-dimensional triangle in $\mathbb
R^2$ and $Q$ is an edge of $P$. In the quotient by
$\operatorname{span}(Q)$, the vector space parallel to $Q$, the images
of the integral points are shown with dots and the image
$P/\operatorname{span}(Q)$ is a line segment of length 3 with $p_Q$ as one of
its endpoints. The integer vector $\mathbf w$ is a minimal integer normal vector
to $p_q$ in $P$, and in addition, we've illustrated two of its preimages as
$\tilde{\mathbf w}_1$ and $\tilde{\mathbf w}_2$. Any affine linear function
which is constant along $Q$ will have the same slopes along $\tilde{\mathbf
w}_1$ and $\tilde{\mathbf w}_2$.}\label{f:slope}
\end{figure}

\begin{const}\label{con:slope}
Let $P$ be a $k$-dimensional polyhedron in $\RR^n$, with rational slopes (or in
a simplex of a weak tropical complex, via the identification of that simplex
with a unimodular simplex). Suppose that $Q$ is a $(k-1)$-dimensional face
of~$P$. We define the \defi{linear span of~$Q$} to be:
\begin{equation*}
\operatorname{span}(Q) = \{c(x-y) \mid x, y \in Q, c \in \RR\},
\end{equation*}
which is a vector subspace of $\RR^n$. In the quotient vector space
$\RR^n/\operatorname{span}(Q)$, the image of $Q$ is a point $p_Q$ and the image
of $P$ is a line segment containing~$p_Q$. In addition, because $P$ has rational
slopes, $\operatorname{span}(Q) \cap \ZZ^n$ is a subgroup of rank $k-1$, so the
image of $\ZZ^n$ in $\RR^n/\operatorname{span}(Q)$ is a discrete lattice.
Therefore, there exists a unique minimal vector $\mathbf w \in \RR^n
/\operatorname{span}(Q)$ which is contained in the image of $\ZZ^n$ and which is
parallel to the image $P/\operatorname{span}(Q)$ and pointing from $p_Q$ into
$P/\operatorname{span}(Q)$,
meaning that $p_Q + \epsilon \mathbf w$ is in $P/\operatorname{span}(Q)$ for
sufficiently small $\epsilon >
0$. (This quotient construction is illustrated in
Figure~\ref{f:slope} for $n = k =2$.)

Now suppose that $\phi$ is an affine linear function on $\RR^n$ which is
constant on $Q$ and has integral slopes. Then, $\phi$ is the pullback of an
affine linear function $\phi'$ on $\RR^n / \operatorname{span}(Q)$. We define
the \defi{slope of $\phi$ on $P$ from $Q$} to be $\phi'(x_Q + \mathbf w) -
\phi'(x_Q)$, which is an integer because $\phi$ has integer slopes and $\mathbf
w$ is in the image of an integer point.
\end{const}

It is important to note that the vector space quotient in
Construction~\ref{con:slope} works only in a single simplex. Therefore, unlike
the case of tropical varieties in $\RR^n$, we cannot define the balancing
condition for formal sums of polyhedra on a weak tropical complex directly in
terms of that quotient. Instead, only the slopes produced by
Construction~\ref{con:slope} are used to define the multiplicity of a PL
function, and from that, the balancing condition.

\begin{const}\label{con:multiplicity}
Let $Z$ be a formal sum of $k$-dimensional polyhedra in an open subset $U$ of a
weak tropical complex $\Delta$. Let $Q \subset U$ be a $(k-1)$-dimensional
polyhedron, and let $P_1, \ldots, P_r$ be the terms of $Z$ which intersect~$Q$
and assume
that $Q$ is a face of each $P_i$. Let $m_1, \ldots, m_r$ be the
coefficients of $P_1, \ldots, P_r$, respectively, in the formal sum $Z$. In
addition, suppose that $\phi$ is a PL function on $U$ such that $\phi$ is
constant on $Q$ and each $P_i$ is contained in a single closed domain of
linearity of $\phi$. Under these assumptions, we define the \defi{multiplicity
of $\phi$ along $Q$ of
$Z$} to be:
\begin{equation*}
\mult_{Q,Z}(\phi) = \sum_{i=1}^r a_i m_i,
\end{equation*}
where $a_i$ is the slope of $\phi$ on $P$ from $Q$, as in
Construction~\ref{con:slope}.
\end{const}

The requirement in Construction~\ref{con:multiplicity} that $Q$ is a face of all
the terms of~$Z$ which it intersects is a strong one, but can be achieved by
subdividing. In particular, for any PL function $\phi$ which is constant on a
$(k-1)$-dimensional polyhedron $Q$, and any formal sum of $k$-dimensional
polyhedra $Z$, we can find a subdivision $Z'$ and a $(k-1)$-dimensional polytope
$Q' \subset Q$ such that any term of $Z'$ which intersects $Q'$ has $Q'$ as a
face. Such subdivisions are necessary when the multiplicity varies across $Q$
because terms of $Z$ may intersect only parts of $Q$ or because the slopes of
$\phi$ may also vary across points of $Q$.

On the other hand, once the requirements in Construction~\ref{con:multiplicity}
for defining the multiplicity are satisfied, further
subdivisions will inherit the same multiplicities, in the following sense.
Suppose that $Z$ is a formal sum of
$k$-dimensional polyhedra in an open subset $U$ of a weak tropical complex. Let
$Q$ be a $(k-1)$-dimensional polyhedron and $\phi$ a PL function such that $Q$
is a face of any term of $Z$ that it intersects, any term of $Z$ is contained in
a domain of linearity of~$\phi$, and $\phi$ is constant on~$Q$, meaning that
Construction~\ref{con:multiplicity} defines a multiplicity $\mult_{Q,Z}(\phi)$.
Then, for any subdivision of $Z'$ of~$Z$, and $Q' \subset Q$, which satisfy the
same conditions so that $\mult_{Q', Z'}(\phi)$ is defined, then
$\mult_{Q',Z'}(\phi)$ will equal $\mult_{Q, Z}(\phi)$. This is because in the
subdivision process, each term $P$ of $Z$ containing $Q$ will contribute exactly
one term $P'$ of $Z'$ containing $Q'$, on which $\phi$ will have the same slope
as it did on~$P$.

\begin{defn}\label{d:balanced}
Let $Z$ be a formal sum of $k$-dimensional polyhedra. Suppose that the
intersection of any two terms of $Z$ is either a face of each or has dimension
at most $k-2$. We say that $Z$ is \defi{balanced} if for each
$(k-1)$-dimensional face~$Q$ of a term of $Z$, and any linear function $\phi$ on
a neighborhood of the interior of $Q$, the multiplicity $\mult_{Q,Z}(\phi)$ is
zero.
\end{defn}

\begin{lem}\label{l:balanced-equiv}
Suppose that $Z$ and $Z'$ are equivalent formal sums of $k$-dimensional
polyhedra, such that the intersection of any two terms in $Z$ is either a face
of each, or has dimension at most $k-2$, and $Z'$ has the same property. Then
$Z$ is balanced if and only if $Z'$ is balanced.
\end{lem}

\begin{proof}
Given formal sums $Z$ and $Z'$ as in the lemma statement, we can find a common
refinement of all their terms, and refine further such that the intersection of
two terms is either a face of each or has dimension at most $k-2$. Thus, it
suffices to prove the theorem if we assume that $Z'$ is obtained by repeatedly
refining the terms of $Z$.

In this case, for any $(k-1)$-dimensional face~$Q$ of $Z$, there is a
$(k-1)$-dimensional face of $Z'$ contained in $Q$, and so if $Z'$ is balanced,
then $Z$ is balanced because of the compatibility between multiplicities and
subdivisions, as noted above.

Conversely, assume that $Z$ is balanced. Given a $(k-1)$-dimensional face $Q'$
of a $Z'$, there are two cases. First, suppose that $Q'$ is contained in a
$(k-1)$-dimensional face $Q$ of $Z$. Then, by the same compatibility of
subdivisions and multiplicities, it is sufficient to show that any linear
function $\phi'$ on a neighborhood of the interior of $Q'$ extends to a linear
function on a neighborhood of the interior of $Q$. By the definition of linear
functions, $\phi'$ extends to a function $\phi$ on the interiors of all
simplices containing $Q'$, which is a neighborhood of the interior of $Q$.
Moreover, $\phi$ is constant on $Q$, so the multiplicity of $\phi$ along $Q$ of
$Z$ is $0$, and thus the multiplicity of $\phi$ along $Q'$ of $Z'$ is also $0$,
which is what we need to show.

The second case is that $Q'$ is contained in a term~$P$ of $Z$, but not as a
face. In this case, $Q'$ is contained in exactly two terms, $P_1'$ and $P_2'$,
of $Z'$, on opposite sides of $Q'$ and both formed by subdividing $P$ and thus
with the same coefficient. In the vector space quotient in
Construction~\ref{con:multiplicity}, $P_1'$ and $P_2'$ will be on opposite sides
of the point $p_{Q'}$, and so have opposite minimal vectors, and therefore,
their contributions to the multiplicity of the linear function $\phi$ along $Q'$
will cancel. We've shown that the multiplicities of linear functions along
$(k-1)$-dimensional faces of terms of $Z'$ all, and so $Z'$ is balanced.
\end{proof}

\begin{defn}\label{d:cycle}
A \defi{$k$-cycle} on $\Delta$ is an equivalence class of formal sums of
$k$-dimensional polyhedra such that any representative $Z$ (or, equivalently,
some representative) where the intersection of two of terms of $Z$ is either a
face of each or has dimension at most $n-2$ is balanced. An $(n-1)$-cycle is
called a \defi{Weil divisor} and a $1$-cycle is called a \defi{curve}.
\end{defn}

\begin{rmk}\label{r:yu-balancing}
A different balancing condition for 1-dimensional cycles on the dual complex of
a degeneration is given in~\cite{yu-balancing}. Yu's balancing is phrased in
terms of numerical classes of curves, but is equivalent to using the formalism
in Definition~\ref{d:cycle}, except with the more restrictive definition of
linear functions given in \cite{yu-compactness}*{Def.~3.2} and described in
Remark~\ref{r:yu-linear}. Since Yu allows fewer linear functions, any curve in
our sense is balanced in Yu's sense, but not necessarily vice versa.
Yu shows that the specialization of any algebraic (even
analytic) curve is balanced in his sense. In contrast, in order for the
specialization of an algebraic curve to be a curve on the weak tropical complex,
we will require additional hypotheses, given in
Corollary~\ref{c:specialization-curve} below.
\end{rmk}

\begin{rmk}\label{r:tropical-balancing}
The balancing condition for curves in tropical geometry is usually expressed as
the vanishing of the sum of the minimal outgoing vectors around a
vertex~$v$~\cite{maclagan-sturmfels}*{Def. 3.3.1}. The vanishing of each component of
this sum is equivalent to the vanishing of the multiplicity of the coordinate
function at~$p$, as in Definition~\ref{d:cycle}. The balancing for
$n$-dimensional varieties for $n \geq 2$ is expressed by first taking the
quotient by the linear span of a $(n-1)$-dimensional simplex, which is
equivalent to considering linear functions which are constant on that simplex.
\end{rmk}

We now define a specialization map from cycles on the algebraic
variety~$\X_\eta$ to the weak tropical complex $\Delta$. Let $Z$ be an
irreducible subvariety of~$\X_\eta$ of dimension $k$ and let $\overline Z$ be
its closure in $\X$. We define the specialization as a formal sum of
$k$-dimensional simplices of~$\Delta$:
\begin{equation}\label{eq:specialization}
\rho(Z) = \sum_{s \in \Delta_k} (\deg C_s \cdot \overline Z ) [s]
= \sum_{s \in \Delta_k}
(\deg (C_{s_0} \cdots C_{s_{k}} \cdot \overline Z)\vert_{C_s})[s].
\end{equation}
Since $\X$ is not a variety over a field, we explain the intersection product
$C_s \cdot \overline Z$ in~(\ref{eq:specialization}). We first take the
intersection intersection of the cycle $\overline Z$ with the succession of
Cartier divisors $C_{s_i}$ as $s_i$ runs over the vertices of the simplex~$s$,
which gives a linear equivalence class of 0-dimensional cycles on the
set-theoretic intersection $C_{s_0} \cap \cdots \cap C_{s_{k}} \cap \overline
Z$. Second, we take the total degree of the
cycles contained in $C_s$, which is one connected component of the intersection
$C_{s_0} \cap \cdots \cap C_{s_{k}}$. We extend the specialization map $\rho$ to
$k$-cycles by linearity.

\begin{defn}\label{d:simplicial}
A \defi{simplicial PL function} is a PL function on~$\Delta$, which takes
integer values at the vertices of $\Delta$ and  whose restriction
to each simplex is linear function, under the identification of a simplex with a
standard unit simplex. In other words, a \defi{simplicial PL function} $\phi$ is
defined by choosing integer values $\phi(v) \in \ZZ$ for each vertex $v \in
\Delta_0$ and interpolating linearly on each simplex. A simplicial PL function
on an open subset $U \subset \Delta$ is any restriction of a simplicial PL
function on~$\Delta$.
\end{defn}

We need the following lemma to relate our definition of linearity to the
algebraic geometry of a degeneration.

\begin{lem}\label{l:linear-degree-zero}
Let $\phi$ be a simplicial PL function on an open set $U \subset \Delta$, and
suppose that $\phi$
is constant on a $k$-dimensional simplex $s$. Let $D_\phi$ be the divisor on
$C_s$ defined as the sum:
\begin{equation*}
D_\phi = \sum_{s' \in \Delta_{k+1}, s' \supset s}
a_{s'}[C_{s'}],
\end{equation*}
where $a_{s'}$ is the slope of $\phi$ along $s'$. Then $\phi$ is linear in a
neighborhood of~$s$ if and only if the intersection $D_\phi \cdot C_r$ in $C_s$
has degree zero for every ridge $r$ containing $s$.
\end{lem}

\begin{proof}
Let $r$ be a ridge of $\Delta$ containing $s$. Then, the pullback $\phi \circ
\pi_r \colon N_r^o \rightarrow \RR$ is the restriction of a simplicial PL
function~$\tilde \phi$ on $N_r$ defined by $\tilde\phi(v) = \phi(s)$ for any
vertex~$v$ such that $\pi_r(v)$ is in $s$, and $\tilde\phi(v) =
\phi(s) + a_{s'}$, where $s'$ denotes the unique simplex containing $s$ and
$\pi_r(v)$, and $a_{s'}$ is the slope of $\phi$ along $s'$.

We now label the vertices of $r$ as $v_1, \ldots v_n$, and label the other
vertices of $N_r$, which are in bijection to the vertices of the link of $r$, as
$w_1, \ldots, w_d$ as in Construction~\ref{con:linear-function}. We can assume
that $v_1, \ldots, v_{k+1}$ are the vertices of $s$. For $i$ from
$1$ to $d$, we let $\tilde s_i'$ denote the $(k+1)$-dimensional simplex of $N_r$
which contains $s$ and $w_i$, and so $\pi_r(\tilde s_1'), \ldots, \pi_r(\tilde
s_d')$ are the $(k+1)$-dimensional simplices of $\Delta$ which contain $s$, but
are not contained in $r$, although possibly with repetition. Thus, we let $s_1',
\ldots, s_l'$
denote the $l \leq d$ distinct $(k+1)$-dimensional simplices containing $s$ but
not contained
in~$r$. Finally, $t_{k+2}, \ldots t_{n}$ denote the $(k+1)$-dimensional
faces of $r$ containing $s$ and the vertices $v_{k+2}, \ldots, v_{n}$,
respectively.

Recall that in
Construction~\ref{con:linear-function}, we mapped the vertices $v_1, \ldots,
v_n$ and $w_1, \ldots, w_d$ to a basis of $\RR^{d+n}$ and therefore $\tilde \phi$ is
the restriction of a linear function on $\RR^{d+n}$. However, it is the pullback
of a linear function on $L_r = \RR^{d+n}/\RR$ if and only if it vanishes on the
line
spanned by the vector $(1, \ldots, 1, -\alpha(v_1, r), \ldots, -\alpha(v_n,r))$,
which means that:
\begin{align*}
0 &= \tilde\phi(w_1) + \cdots + \tilde\phi(w_d) - \alpha(v_1, r) \tilde\phi(v_1) -
\cdots - \alpha(v_n, r) \tilde \phi(v_n) \\
&= \phi(s) + a_{\pi(\tilde s_1')} + \cdots + \phi(s) + a_{\pi(\tilde s_d')}
- \alpha(v_1,r) \phi(s) - \cdots - \alpha(v_{k+1},r) \phi(s) \\
&\qquad\qquad - \alpha(v_{k+2},r)
(\phi(s) + a_{t_{k+2}}) - \cdots - \alpha(v_{n},r) (\phi(s) + a_{t_n}) \\
&= a_{\pi_r(s_1')} + \cdots + a_{\pi_r(s_d')} - \alpha(v_{k+2}, r) a_{t_{k+2}} - \cdots -
\alpha(v_n, r) a_{t_{n}} \\
\intertext{by the relation (\ref{eq:constraint}) in the definition of a weak tropical
complex}
&= a_{s_1'} (\# C_{s_1'} \cap C_r) + \cdots + a_{s_l'} (\# C_{s_l'} \cap C_r)
+ a_{t_{k+2}} \deg C_{v_{k+2}} \cdot C_r + \\
&\qquad\qquad\qquad \cdots + a_{t_{n}} \deg C_{v_n} \cdot C_r, \\
\intertext{where the intersections in the last set of terms are taken in $\X$}
&= a_{s_1'} \deg C_{s_1'} \cdot C_r  + \cdots + a_{s_l'} \deg C_{s_l'} \cdot C_r
+ a_{t_{k+2}} \deg C_{t_{k+2}} \cdot C_r + \\
&\qquad\qquad\qquad \cdots + a_{t_{n}} \deg C_{t_{n}} \cdot C_r, \\
\intertext{where the intersections are now taken inside $C_s$, because
$C_{s_i'}$ and $C_r$ intersect transversely at a finite number of points, and
because $C_{v_i} \cap C_s = C_{s_{d+i}'}$ and that intersection is transverse}
&= \deg D_\phi \cdot C_r
\end{align*}
by linearity and the definition of $D_\phi$.

Therefore, $\phi$ is linear on a neighborhood of $r$ if and only if $\deg D_\phi
\cdot C_r = 0$. By definition, $\phi$ is linear on a neighborhood of $s$ if and
only if it is linear on a neighborhood of each ridge containing $s$, and so
we've proved the lemma statement.
\end{proof}

The specialization defined in (\ref{eq:specialization}) is not necessarily
balanced without additional hypotheses on the degeneration. The following lemma
gives a somewhat technical criterion for a specialization to be a $k$-cycle,
which we will subsequently specialize to the case of curves, which are a key
case of interest in this paper.

\begin{lem}\label{l:cycle-criterion}
Let $Z$ be a $k$-dimensional subvariety of $\X_\eta$. Then $\rho(Z)$ is a
$k$-cycle if and only there exist rational numbers $t_{s,r}$,  indexed by  the
pair of a $(k-1)$-dimensional simplex $s$ in $\Delta$ and a ridge $r$ contianing
$s$, such that, for any $k$-dimensional simplex $s'$ containing $s$:
\begin{equation*}
\deg C_{s'} \cdot (C_s \cdot \overline Z) =
\sum_{r \in \Delta_{n-1}, r \supset s}^m t_{s,r} \deg C_{s'} \cdot
C_{r}.
\end{equation*}
Here, the sum is over all ridges $r$ containing $s$ and the intersection $C_s
\cdot \overline Z$ is taken in $\X$, but all other intersection products are
taken inside $C_s$.
\end{lem}

\begin{proof}
First, suppose that $s$ is a $(k-1)$-dimensional simplex such that there do not
exist rational numbers $t_{s,r}$ as in the lemma statement. Then, by linear
algebra, there exists a rational linear combination $D = a_1 C_{s_1'} + \cdots +
a_l C_{s_l'}$, such that $D \cdot C_r = 0$ for all ridges $r$ containing $s$, but
$D \cdot (C_s \cdot \overline Z) \neq 0$. By scaling, we can assume that the
$a_i$ are integers, and then construct the function $\phi$ on an open
neighborhood of the interior of $s$ which is identically zero on $s$ and has
slope $a_i$ on $s_i'$. Therefore, $D = D_\phi$, in the notation of
Lemma~\ref{l:linear-degree-zero}, and thus $\phi$ is a linear function by that
lemma.

In addition, we want to show that the multiplicity of $\phi$ at $s$ is non-zero.
Let $m_1, \ldots, m_l$ denote the multiplicities of $\rho(Z)$ at $s_1', \ldots,
s_l'$, respectively, which denote the $k$-dimensional simplices containing $s$.
Then,
\begin{equation}\label{eq:multiplicity-degree}
\begin{split}
\mult_{\rho(Z),s}(\phi) &= \sum a_i m_i
= \sum a_i \deg C_{s_i'} \cdot \overline Z \\
&= \deg \left( \sum a_i C_{s_i'} \right) \cdot (C_s \cdot \overline Z) 
= \deg D_\phi \cdot (C_s \cdot \overline Z), 
\end{split}
\end{equation}
which is non-zero by assumption.

Conversely, suppose that there exist rational numbers $t_{s,r}$, as in the lemma
statement.
We will show that
$\rho(Z)$ is balanced along a given $(k-1)$-dimensional simplex $s$. As before,
$s_1', \ldots, s_l'$ are the $k$-dimensional simplices containing $s$, and $m_1,
\ldots, m_l$ are their
respective multiplicities in $\rho(Z)$. Let $\phi$ be a linear function on an
open set which meets $s$. Then $\phi$ is a simplicial PL function, and we let
$a_i$ denote the slope of $\phi$ along $s_i'$. Let $D_\phi$ denote the sum
divisor $\sum_i a_i [C_{s_i'}]$, as in Lemma~\ref{l:linear-degree-zero}. The
computation in~(\ref{eq:multiplicity-degree}) applies and so we have:
\begin{align*}
\mult_{s, \rho(Z)}(\phi)
&= \deg D_\phi \cdot (C_s \cdot \overline Z) \\
&= \deg D_\phi \cdot \sum_{r \in \Delta_{n-1}, r \supset s} t_{s,r} [C_{r}], \\
\intertext{by our assumption on the $t_{s,r}$ and because $D_\phi$ is a linear
combination of the~$C_{s_i'}$.}
&= \sum_{r \in \Delta_{n-1}, r \supset s} t_{s,r} \deg D_\phi \cdot C_{r} = 0
\end{align*}
by Lemma~\ref{l:linear-degree-zero} and because $\phi$ is linear. Therefore,
$\rho(Z)$ is balanced.
\end{proof}

\begin{cor}\label{c:specialization-curve}
If $\X$ is a numerically faithful degeneration, and $Z$ is a $1$-cycle in
$\X_\eta$, then $\rho(Z)$ is a curve.
\end{cor}

\begin{proof}
Let $v$ be a vertex of the tropical complex. By the definition of numerically
faithful, the curves $C_r$, as $r$ ranges over ridges containing $v$ generate
the vector space of numerical classes of curves, meaning that a rational linear
combination of the $C_r$ is numerically equivalent to the curve $C_v \cap
\overline Z$. Such a linear combination satisfies the criterion of
Lemma~\ref{l:cycle-criterion}, so $\rho(Z)$ is a curve.
\end{proof}

In this paper, we're interested in specializations of curves and divisors, but
we note that an analogue of Corollary~\ref{c:specialization-curve} is true for
specializations of $k$-cycles if we instead assume that the $1$-dimensional
strata generate the vector space of numerical curve classes in the stratum $C_s$
for each $(k-1)$-simplex in~$\Delta$. In particular, for the specialization of a
divisor on~$\X_\eta$ to be a Weil divisor, the criterion would be that the vector
space of numerical classes of curves on each surface $C_q$ is generated by the
1-strata in $C_q$. However, the weaker condition of robustness in dimension~2 is
also sufficient:

\begin{prop}\label{p:specialization-divisor}
If $\X$ is robust in dimension~$2$, and $Z$ is a divisor on~$\X_\eta$, then
$\rho(Z)$ is a Weil divisor on~$\Delta$.
\end{prop}

\begin{proof}
Suppose that $\X$ is robust in dimension~$2$, and let $Z$ be a divisor
on~$X_\eta$. Fix a $(n-2)$-dimensional simplex~$q$ of~$\Delta$. By
Proposition~\ref{p:tropical-complex-robustness}, there exists a divisor $B$ of
the form $B = \sum_r b_r [C_r]$ on $C_q$, where the sum is over ridges $r$
containing $q$, such that $\deg B^2 > 0$.

In addition, suppose that $\phi$ is a linear function on a neighborhood of the
interior of $q$. Let $D_\phi = \sum_r a_r [C_r]$ be the divisor on $C_q$, where,
as in Lemma~\ref{l:linear-degree-zero}, the sum is over ridges $r$ containing
$q$ and $a_r$ is the slope of $\phi$ along $r$. Then, by
Lemma~\ref{l:linear-degree-zero}, $D_\phi \cdot C_r$ has degree zero for any
ridge $r$ containing $q$, and so $\deg D_\phi \cdot B = \deg D_\phi^2 = 0$. If
we consider the three divisors, $B$, $D_\phi$, and $\overline Z \cdot C_q$ on
$C_q$, then the matrix of the intersection pairing of these divisors, in that
order, on $C_q$
is:
\begin{equation*}
\begin{bmatrix}
> 0 & 0 & * \\
0 & 0 & c \\
* & c & * 
\end{bmatrix}.
\end{equation*}
If $c = \deg (\overline Z \cdot C_q) \cdot D_\phi$ is non-zero, then this matrix
has two positive eigenvalues, no matter the values in the entries marked $*$.
Therefore, by the Hodge index theorem, $c = 0$.

Finally, we expand the expression  for $D_\phi$,
\begin{align*}
c &= \deg (\overline Z \cdot C_q) \cdot D_\phi
= \sum_{r \in \Delta_{n-1}, q \subset r} a_r \deg (\overline Z \cdot C_q) \cdot
C_r  \\
&= \sum_{r \in \Delta_{n-1}, q \subset r} a_r \deg \overline Z \cdot C_r, \\
\intertext{where the intersection product is now taken inside $\X$ instead of $C_q$,}
&= \mult_{q, \rho(Z)}(\phi),
\end{align*}
by the definition of the multiplicity and because $\overline Z \cdot C_r$ is
the coefficient of $r$ in $\rho(Z)$.
\end{proof}

\section{Cartier divisors}\label{s:cartier}

In this section, we introduce a formal sum of polyhedra to any PL function,
which is the combinatorial analogue of the divisor associated
to a rational function on an algebraic variety. Building on this definition, we
define Cartier divisors as formal sums of polyhedra which are locally defined by
PL functions.

Construction~\ref{con:multiplicity} assigned a multiplicity to each boundary
between domains of linearity of a PL function, but under the assumption that
the function was constant along that boundary. In order to assign multiplicities
to arbitrary PL functions on open subsets of~$\Delta$, the key additional tool
is to add the assumption that the multiplicity of a linear function is 0, which
appears as (iii) in the following proposition. The normalization in (iv) takes
the place of the slopes used in Construction~\ref{con:multiplicity}.

\begin{prop}\label{p:divisor-pl}
For any $n$-dimensional weak tropical complex~$\Delta$, there is a unique
function, denoted
$\div$, which takes a PL function $\phi$ on any open set $U \subset \Delta$ to
an equivalence class of formal sum of polyhedra in the same subset $U$, with the
following properties:
\begin{enumerate}
\item[(i)] For any PL functions $\phi$ and~$\phi'$ on the same open set $U$,
$\div(\phi + \phi') = \div(\phi) + \div(\phi')$.
\item[(ii)] If $V \subset U$ is open and $\phi$ is a PL function on $U$,
then $\div(\phi \vert_V)
= \div(\phi)\vert_V$.
\item[(iii)] The function~$\phi$ is linear 
 if and only if $\div(\phi)$ is trivial.
\item[(iv)] Suppose that $\phi$ is identically zero outside of a single facet~$f$
of~$\Delta$, on which it is defined by:
\begin{equation}\label{eq:piecewise}
\phi(x) = \max \{ \lambda \cdot x, 0\},
\end{equation}
where $x$ is a coordinate vector in standard unit simplex identified with~$f$,
and $\lambda$ is an integral vector whose entries have no
non-trivial common divisor. Then, $\div(\phi) = [\{x \in f \cap U \mid \lambda
\cdot x = 0\}]$.
\item[(v)]\label{it:ridge} If $r$ is a ridge not contained in any facet, $U$ is
the interior of~$r$, and $\phi$ is a simplicial PL function on~$r$, then
\begin{equation*}
\div(\phi) = -\bigg(\sum_{v \in \tilde r} \phi(v) \alpha(v, r)\bigg)[r],
\end{equation*}
where we extend $\phi$ by continuity to a linear function on the
$(n-1)$-dimensional simplex $\tilde r$, which is possible because $\phi$ is a
simplicial PL function.
\end{enumerate}
\end{prop}

\begin{proof}
By property (ii), we can work locally. Since a PL function, restricted to the
interior of a domain of linearity, is linear in the sense of
Construction~\ref{con:linear-function}, it suffices to check that $\div(\phi)$ is
uniquely defined in two cases: along the boundary between domains of linearity
within a facet and along simplices of dimension at most $(n-1)$. Moreover, since
$\div(\phi)$ is a formal sum of $(n-1)$-dimensional polyhedra, we can ignore
sets of dimension $n-2$ or less, and we can refine our two cases to the boundary
between two domains of linearity in a facet and along ridges of~$\Delta$.

In the second case, when, in addition, the ridge~$r$ is not contained in any
facets, we can fix a domain of linearity~$U$ in~$r$. Then we can extend
$\phi\vert_U$ to a unique simplicial PL function on~$\tilde r$.
Then, the divisor of this extension is determined by property~(v), which gives
us the divisor of~$\phi\vert_U$ by property~(ii). Note that, since the slopes of
$\phi$ are required to be integral, and the sum of the $\alpha(v,r)$ is $0$, the
coefficient in property~(v) is integral. Also, that formula is clearly linear
in~$\phi$, and it gives $0$ if and only if $\phi$ is linear, since
Construction~\ref{con:linear-function} involves the quotient by the vector
consisting of the $\alpha(v,r)$.

We now address the remaining cases, that of the boundary between two domains of
linearity in a single facet and of a ridge~$r$ contained in one or more facet.
In the first case, we define $\phi'$ to be the linear extension of $\phi$ from
one of the domains of linearity. In the second case, we can also find a linear
function~$\phi'$ which locally agrees with $\phi$ on all but one of the facets
containing~$r$. The reason is that the map~$\phi_r$ from
Construction~\ref{con:linear-function} imposes a single relation so that if we
remove one facet, the images of the remaining vertices are affinely independent.
In either case, $\phi - \phi'$ is supported on a single facet, and is linear on
its support. Thus, it can be written as an integral multiple of function as
in~(\ref{eq:piecewise}), so properties~(v) and then~(i) compute the divisor.

In both of these cases, the construction of $\phi'$ is linear in $\phi$, and so
the resulting multiplicity according to (iv) is linear in $\phi$ as well, so
long as the we choose the same support for $\phi'$. Thus, the only remaining
point is to check that the multiplicity computation is independent of the choice
of $\phi'$.

In other words, we need to
show that, if two distinct functions $\phi$ and~$\phi'$ both have the same form
as the function~(\ref{eq:piecewise}), and $\phi - \phi'$ is linear, then both
functions define the same divisor. Thus, the linear function $\phi - \phi'$ must
agree with $\phi$ where it is non-zero. If $\phi$ is defined on the interior of
a facet, then this means that $\phi - \phi'$ must be $\lambda \cdot x$, where
$\lambda$ is the slope vector of the non-zero part of~$\phi$, as
in~(\ref{eq:piecewise}), so $\phi' = \max\{-\lambda \cdot x, 0\}$, which also
defines a divisor of multiplicity~$1$ along the set where $\lambda \cdot x = 0$.
On the other hand, if $\phi$ is defining a divisor contained in a ridge~$r$,
then $\phi$ and~$\phi'$ are non-zero on all but one of the facets
containing~$r$, and differ in which facet that is. However, since the vector
defining the quotient in Construction~\ref{con:linear-function} has an entry
of~$1$ for each facet containing~$r$, the $\phi$ and $\phi'$ both have
slope~$1$, so they define the same divisor when applying property~(iv).
\end{proof}

\begin{defn}\label{d:cartier}
On an $n$-dimensional weak tropical complex~$\Delta$, a
\defi{Cartier divisor} $D$ is an equivalence class of formal sums of
$(n-1)$-dimensional polyhedra such that $D$
is locally defined by a PL function, in the sense that there exists an open
cover $U_1, \ldots, U_m$ of $\Delta$ and PL functions $\phi_i$ on $U_i$ such
that $D \vert_{U_i} = \div(\phi_i)$ for each $i$.

A \defi{$\QQ$-Cartier divisor} is an equivalence class of formal (integral) sums
of $(n-1)$-dimensional polyhedra such that $mD$ is a Cartier divisor for some
positive integer $m$.
\end{defn}

\begin{prop}\label{p:cartier-weil}
If $\Delta$ is a 2-dimensional weak tropical complex, then a formal sum of
1-dimensional polyhedra is a Weil divisor if and only if it is a $\QQ$-Cartier
divisor.
\end{prop}

\begin{proof}
See \cite{cartwright-surfaces}*{Prop. 2.8}.
\end{proof}

\begin{rmk}
The proof of Proposition~\ref{p:cartier-weil} given
in~\cite{cartwright-surfaces} can be generalized to show that Weil divisors on a
$n$-dimensional weak tropical complex are $\QQ$-Cartier in codimension 1,
i.e.\ they are $\QQ$-Cartier except on a set of dimension at most $n-3$, but, the
converse does not hold in general. In earlier versions of this paper,
``$\QQ$-Cartier in codimension 1'' was given as the definition of a Weil
divisor, but Definition~\ref{d:cycle} is more natural from the perspective of
tropical geometry, and from how the condition is used in the proofs.
\end{rmk}

\begin{figure}
\includegraphics{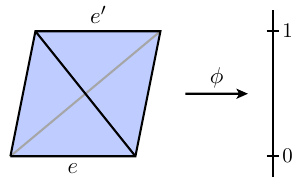}
\caption{A PL function~$\phi$ on a tropical complex consisting of a tetrahedron
with $\alpha(v,r) = 1$ for every vertex~$v$ in every edge~$r$. The divisor
associated to~$\phi$ is $2[e] - 2[e']$.}
\label{f:tetrahedron_pl}
\end{figure}

\begin{ex}\label{ex:tetrahedron-div}
We again consider the tetrahedron~$\Delta$ with all structure constants
$\alpha(v,r)$ equal to
$1$, as in Example~\ref{ex:k3-robust}. Consider the PL function~$\phi$ depicted
in Figure~\ref{f:tetrahedron_pl} and we explain the computation of $\div(\phi)$
as in Proposition~\ref{p:divisor-pl}. Near the edge~$e$, $\phi$ is the sum of a
function increasing along one of the facets containing $e$ and constant on the
other, and a function with the opposite behavior. By (iv) from
Proposition~\ref{p:divisor-pl}, each of these functions has divisor $[e]$, which
shows that $\div(\phi)$ has a term $2[e]$. After subtracting the constant $1$
from $\phi$, a similar computation gives a coefficient of $-2$ for
the edge $[e']$. Finally, one can check that $\phi$ is linear along each of the
other 4 ridges of $\Delta$, showing that $\div(\phi) = 2[e] - 2[e']$.

From this computation, it immediately follows that $2[e]$ is a Cartier divisor
and that $[e]$ is a $\QQ$-Cartier divisor. In Example~\ref{ex:k3-intersection},
we will see that $[e]$ is not a Cartier divisor. By
Proposition~\ref{p:cartier-weil}, both of these are Weil divisors, which can
also be checked directly from the fact that the only linear functions in a
neighborhood of the vertices of~$\Delta$ are the constant functions.
\end{ex}

The following construction gives an intersection product between Cartier
divisors and curves. Note that for formal sums of 0-dimensional polyhedra, the
equivalence relation and balancing condition are both trivial, so a 0-cycle is
the same thing as a formal sum of points.

\begin{const}\label{con:intersection}
Let $D$ be a Cartier divisor and $C$ a curve in a weak tropical
complex~$\Delta$. By definition,
there exists an open cover $\{U_i\}$ of $\Delta$ such that
on each~$U_i$, the divisor $D$ is defined by a PL function $\phi_i$. Then, we
define the \defi{intersection product $D \cdot C$} to be the 0-cycle such that,
in each $U_i$, we have the restriction:
\begin{equation}\label{eq:intersection-const}
(D \cdot C)\vert_{U_i} = \sum_{p \in U_i \cap C} \mult_{p,C}(\phi_i) [p],
\end{equation}
where, in order to use the multiplicity from
Construction~\ref{con:multiplicity}, we've implicitly chosen a representative
formal sum of 1-dimensional polyhedra from the equivalence class $C$ such that
any edge containing $p$ has $p$ as an endpoint.
The sum of the coefficients in $D \cdot C$ is its degree, denoted
$\deg D \cdot C$.

If $D$ is a $\QQ$-Cartier divisor, and $m$ is a positive integer such that $mD$
is Cartier, then we define $D \cdot C$ to be $\frac{1}{m} (mD) \cdot C$, which
is a formal sum of points with rational coefficients.
\end{const}

\begin{prop}\label{p:intersection}
If $D$ is a Cartier (resp.\ $\QQ$-Cartier) divisor on a weak tropical complex
$\Delta$ and $C$ is a curve, then
intersection $D \cdot C$ from Construction~\ref{con:intersection} is a
well-defined formal sum of points (resp.\ formal rational sum of points).
\end{prop}

\begin{proof}
We first show that the sum given in~(\ref{eq:intersection-const}) is finite.
Let $\phi_i$ be one of the PL functions defining the Cartier divisor~$D$. Then,
at any point~$p$ in the interior of a segment $C$ and in the interior of a
domain of linearity for~$\phi_i$, there will be two outgoing slopes from $p$,
which will cancel. Since $C$ has finitely many segments and $\phi_i$ has
finitely many domains of linearity, this shows that on each of the finitely many
open sets~$U_i$, the sum~(\ref{eq:intersection-const}) has finitely many
non-zero terms, and so $D \cdot C$ is a finite formal sum.

There are four sets of choices made in Construction~\ref{con:intersection}:
first, the open sets~$U_i$, second for the local defining equations, third, the
choice of subdivision in order to define the multiplicity as in
Construction~\ref{con:multiplicity}, and fourth,
the multiple $m$ in the case when $D$ is $\QQ$-Cartier. Because
the divisor of a PL function is a local condition, the $U_i$ can be refined and
the defining equations replaced by their restrictions without changing the
intersection. 

Now suppose we are computing
with the same open cover, but two different sets, $\phi_i$ and $\phi_i'$ of
local defining equations for~$D$. Then their difference $\phi_i - \phi_i'$ has
trivial divisor on each $U_i$, which means that $\phi_i - \phi_i'$
is a linear function by Proposition~\ref{p:divisor-pl}(iii).
Therefore, by the definition of a curve, $\mult_{p,C}((\phi_i - \phi_i')) =
0$ for all points $p \in C$, and thus $\mult_{p,C}(\phi_i) =
\mult_{p,C}(\phi_i')$. Thus, $D \cdot C$ is a well-defined $0$-cycle.

Third, the multiplicity is independent of the choice of a formal sum of
polyhedra in the equivalence class $C$ by the discussion before
Definition~\ref{d:balanced}. Finally, if we choose a larger multiple $km$, we
can compute $(kmD) \cdot C$ using the PL functions $k \phi_i$. Since the
multiplicity scales with the PL function, $(kmD) \cdot C = k (mD) \cdot C$, and
so the computation of $D \cdot C$ is consistent.
\end{proof}

\begin{ex}\label{ex:k3-intersection}
Let $\Delta$ be the tetrahedron with all structure constants equal to~$1$, as in
Example~\ref{ex:k3-robust} and~\ref{ex:tetrahedron-div}. Let $e$ be an edge, and
then $[e]$ is a $\QQ$-Cartier divisor and $2[e]$ is a Cartier divisor, by
Example~\ref{ex:tetrahedron-div}, and thus a curve by
Proposition~\ref{p:cartier-weil}. Let $e''$ be another edge of $\Delta$ which
shares a vertex $v$ in common with $e$. By symmetry, $[e'']$ is also a curve,
and so we can compute $[e] \cdot [e'']$ using
Construction~\ref{con:intersection}.

Let $\phi$ be the PL function from Example~\ref{ex:tetrahedron-div}, so that
$2[e]$ can be defined by $\phi$ in a neighborhood of $e$, and by the constant
$0$ function on the complement of $e$. Then, $\phi \vert_{e''}$ is increasing
with slope $1$ on a neighborhood of $v$, and so $\mult_{v,e''}(\phi) =1$.
Since $\phi$ is linear away from $v$, and the constant zero function is linear
everywhere, this is the only multiplicity and thus
$2[e] \cdot [e''] = [v]$.
Therefore, $[e] \cdot [e''] = (1/2)[v]$, which is not an integral sum, and thus by
Proposition~\ref{p:intersection}, $[e]$ is not a Cartier divisor.
\end{ex}

To relate the divisor of a PL function with algebraic geometry, we have the
following generalization of Lemma~\ref{l:linear-degree-zero} for simplicial PL
functions in a neighborhood of a vertex.

\begin{lem}\label{l:divisor-intersection}
Let $\phi$ be a simplicial PL function on a neighborhood~$U$ of a vertex $v$ in a
weak simplicial complex~$\Delta$. For any edge $e$ containing $v$, let $a_e$
denote the slope of $\phi$ along $e$. As in Lemma~\ref{l:linear-degree-zero}, we
define $D_\phi = \sum a_{e} [C_e]$, where the sum is taken over all edges
containing $v$. Then,
\begin{equation*}
\div(\phi) =
\sum_{r \in \Delta_{n-1}, v \in r} (\deg D_\phi \cdot C_r) [r \cap U]
\end{equation*}
\end{lem}

\begin{proof}
Let $r$ be any ridge containing $v$ and we first suppose that $r$ is contained
in at least one facet, which we denote by $f$. Let $\psi$ be the PL function on a
neighborhood of the interior of $r$ which is identically $0$ on $r$ and all of
the facets containing $r$ except for $f$, on which it is increasing with slope
equal to $\deg D_\phi \cdot C_r$. Then, the divisor $D_\psi$ on $C_r$ associated
to $\phi$ is $\deg D_\phi \cdot C_r$ copies of a single point, and so the
difference
$\phi-\psi$ is a linear function on the neighborhood of the interior of $r$, by
Lemma~\ref{l:linear-degree-zero}. Therefore, the multiplicities of $\phi$ and
$\psi$ along $r$ agree by properties (iii) and (i) from
Proposition~\ref{p:divisor-pl}. However, $\psi$ is $\deg D_\phi \cdot C_r$ times
the function described in (iv) of Proposition~\ref{p:divisor-pl}, and so the
multiplicity of $\psi$, and thus $\phi$, along $r$ is $\deg D_\phi \cdot C_r$,
as desired.

Now suppose that $r$ is not contained in any facet. Then, using (ii) from
Proposition~\ref{p:divisor-pl}, the multiplicity of $r \cap U$ is:
\begin{align*}
-\sum_{w \in \tilde r} \phi(w) \alpha(w,r)
&= -\sum_{w \in \tilde r} (\phi(w) + a_{\pi_r(e_w)}) \alpha(w, r),
\intertext{where $e_w$ denotes the edge of $\tilde r$ between $v$ and $w$,}
&= -\sum_{w \in \tilde r} a_{\pi_r(e_w)} \alpha(w,r)
\intertext{by the relation (\ref{eq:constraint}) in the definition of a weak
tropical complex}
&= \sum_{w \in \tilde r} a_{\pi_r(e_w)} \deg C_{\pi_r(w)} \cdot C_r \\
&= \deg D_\phi \cdot C_r
\end{align*}
by linearity of the intersection product and the fact that the only components
of $D_\phi$ which intersect $r$ are those corresponding to the vertices of $r$.
\end{proof}

We now come to the proof of Theorem~\ref{t:intersect-intro}, which shows that the
intersection product on tropical complexes can be used to compute intersection
numbers on algebraic varieties, under certain conditions.

\begin{thm}\label{t:main-body}
Suppose that $\X$ is a numerically faithful degeneration with tropical
complex~$\Delta$. If $D$ is a divisor on $\X_\eta$, then $\rho(D)$ is a
$\QQ$-Cartier divisor on~$\Delta$. If, in addition, $C$ is a curve on $\X_\eta$,
then $\deg \rho(D) \cdot \rho(C) = \deg D \cdot C$.
\end{thm}

\begin{proof}
It suffices to check that $\rho(D)$ is a $\QQ$-Cartier divisor in
a neighborhood~$U$ of each vertex~$v$ of~$\Delta$. We let $\overline D$ be the
closure of~$D$ in~$\X$. Then $\overline D$ is a divisor on $\X$, whose
restriction $D \cap C_v$ and is a divisor on $C_v$,
numerically equivalent to some rational linear combination $\sum a_e [C_{e}]$
where the edges $e$ contain $v$, by the numerically faithful assumption.
Since intersection theory on both $\X_\eta$ and on $\Delta$ is linear, we can
replace $D$ by an integer multiple such that the $a_e$ are all integers, and it
suffices to prove the theorem in this case.

We let $\phi$ be the simplicial PL function which is $0$ at $v$ and has slope
$a_e$ along each edge $e$. Then, by Lemma~\ref{l:divisor-intersection},
\begin{align*}
\div(\phi)
&= \sum_{r \in \Delta_{n-1}, v \in r} (\deg D_\phi \cdot C_r) [r \cap U] \\
&= \sum_{r \in \Delta_{n-1}, v \in r} (\deg (\overline D \cap C_v) \cdot C_r) [r
\cap U]
= \rho(D) \cap U,
\end{align*}
because $D_\phi$ is numerically equivalent to $\overline D \cap C_v$. Therefore,
$\rho(D)$ is a $\QQ$-Cartier divisor.

Now, we compute $\rho(D) \cdot \rho(C)$, on the same neighborhood $U$, and using
the same PL function $\phi$. If $\rho(C) \cap U$ is the sum $\sum m_i [e_i]$, then
the multiplicity of $\rho(D) \cdot \rho(C)$ at $v$ is:
\begin{align*}
\mult_{v,\rho(C)}(\phi) &= \sum_e a_e m_e 
= \sum_e a_e \deg C_e \cdot \overline C 
= \deg (\overline D \cap C_v) \cdot \overline C,
\end{align*}
by the assumed numerical equivalence. Therefore, $\mult_{v,\rho(C)}(\phi)$
computes the degree of $\overline D \cap \overline C$ on the variety $C_v$. By
\cite{fulton}*{Prop. 20.3(b)}, intersection numbers are preserved under
specialization, so we can recover the degree of $D \cdot C$ by summing the
degrees of $\overline D \cap \overline C$ over the components $C_v$, which is
exactly what $\deg \rho(D) \cdot \rho(C)$ computes.
\end{proof}

\section{Linear and algebraic equivalence}\label{s:equiv}

In this section, we define linear equivalence of divisors on a weak tropical
complex, and define the Picard group to be the group of Cartier divisors up
to linear equivalence. We investigate the structure of the Picard group by
introducing algebraic equivalence, which gives us tropical analogues of the
N\'eron-Severi group and the Jacobian. Each of these groups has an
interpretation in terms of sheaf cohomology.

\begin{defn}
Let $\Delta$ be a weak tropical complex. A \defi{principal divisor} on $\Delta$
is a Cartier divisor which is the divisor $\div(\phi)$ of some PL function~$\phi$. Two Weil
divisors or two Cartier divisors on~$\Delta$ are \defi{linearly equivalent} if
their difference is a principal divisor. The \defi{Picard group} $\Pic(\Delta)$
is the group of Cartier divisors on $\Delta$ modulo principal divisors.
\end{defn}

\begin{prop}\label{p:linear-equivalence-trop}
If $Z$ and $Z'$ are linearly equivalent divisors on the generic fiber $\X_\eta$
of a degeneration $\X$ which is robust in dimension~$2$, then the
specializations $\rho(Z)$ and $\rho(Z')$ are linearly equivalent Weil divisors
on the tropical complex~$\Delta$.
\end{prop}

\begin{proof}
If $Z$ and $Z'$ are linearly equivalent divisors on $\X_\eta$, then there exists
a rational function~$f$ on~$\X_\eta$ whose divisor is $Z - Z'$. Regarding $f$ as
a rational function on $\X$, its divisor is $\overline Z - \overline Z' -
\sum_{f \in \Delta_0} b_v [C_v]$, for some integers $b_v$. Define $\phi$ to be
the simplicial PL function defined by setting $\phi(v) = b_v$ and extending
linearly on each simplex, and we claim that $\rho(Z) - \rho(Z') = \div(\phi)$.

Let $r$ be a ridge of $\Delta$, and the multiplicity of $r$ in $\rho(Z) -
\rho(Z')$ is
\begin{equation*}
\deg\overline Z \cdot C_r - \deg\overline Z' \cdot C_r =
\deg(\overline Z - \overline Z') \cdot C_r =
\deg \left(\sum_{v \in \Delta_0} b_v [C_v]\right) \cdot C_r,
\end{equation*}
because the degree of a principal divisor is zero. In addition, choose $u$ to be
an arbitrary vertex of~$r$. Since $\X_0$ is reduced, $\sum_{v \in \Delta_0}
[C_v]$ is a principal divisor, defined by the uniformizer of $R$, and so
\begin{equation*}
\deg \left(\sum_{v \in \Delta_0} b_v [C_v] \right) \cdot C_r =
\deg \left(\sum_{v \in \Delta_0} (b_v - b_u)[C_v]\right) \cdot C_r.
\end{equation*}
Note that $\sum_{v \in \Delta_0} (b_v - b_u) [C_v]$ does not contain the
component $C_u$ and the other components either intersect $C_u$ transversely or
not at all. Therefore, the restriction of the divisor $\sum_{v \in \Delta_0}
(b_v - b_u) [C_v]$ to the component $C_u$ is $\sum_{e \in \Delta_1, u \in e} a_e
[C_e]$, where $a_e$ is defined as $b_v - b_u$, if $v$ is the endpoint of
$e$ other than~$u$. Since $a_e$ is the slope of $\phi$ along $e$, $\sum_{e \in
\Delta_1, u \in e} a_e [C_e]$ is the same as $D_\phi$ from
Lemma~\ref{l:divisor-intersection}, and, by that lemma, the degree of $D_\phi
\cdot C_r$ is the multiplicity of $r$ in $\div(\phi)$. Therefore, $\rho(Z) -
\rho(Z')$ is the principal divisor defined by $\phi$.
\end{proof}

\begin{defn}\label{d:alg-trivial-equiv}
We say that a Cartier divisor $D$ on a weak tropical complex~$\Delta$ is
\defi{algebraically trivial} if there exists an open cover $U_1, \ldots, U_m
\subset \Delta$ and PL functions~$\phi_i$ on $U_i$ such that $D \vert_{U_i} =
\div(\phi_i)$ and such that the difference $\phi_i\vert_{U_i \cap U_j} - \phi_j
\vert_{U_i \cap U_j}$ is locally constant for each pair $i$ and~$j$.

Two Weil divisors or two Cartier divisors are \defi{algebraically equivalent} if
their difference is algebraically trivial. We define the \defi{N\'eron-Severi}
group $\NS(\Delta)$ to be the group of Cartier divisors modulo algebraic
equivalence, and the \defi{Jacobian} $\Jac(\Delta)$ to be the group of
algebraically trivial divisors modulo linear equivalence.
\end{defn}

It is immediate from the definitions that principal divisors are algebraically
trivial, and so linearly equivalent divisors are algebraically equivalent.
In addition, the sum of algebraically trivial divisors is algebraically trivial,
by taking a common refinement of their open covers.
Therefore, we have a short exact sequence of Abelian groups
\begin{equation}\label{eq:jac-pic-ns}
0 \rightarrow \Jac(\Delta) \rightarrow \Pic(\Delta) \rightarrow \NS(\Delta)
\rightarrow 0.
\end{equation}
When $\Delta$ is 1-dimensional, and thus a tropical curve, the Jacobian is
usually defined as the group of divisor classes of degree $0$, where the degree
is the degree of $D \cdot \Delta$, as in Construction~\ref{con:intersection}.
However, our definition of the Jacobian is compatible with the usual one:

\begin{prop}\label{p:alg-equiv-curve}
If $\Delta$ is a 1-dimensional tropical complex, then a divisor~$D$ is
algebraically trivial if and only if it has degree $0$. Therefore, $\NS(\Delta)
\isom \ZZ$ and $\Jac(\Delta)$ is the group of degree $0$ divisors modulo linear
equivalence.
\end{prop}

\begin{proof}
\begin{figure}
\includegraphics{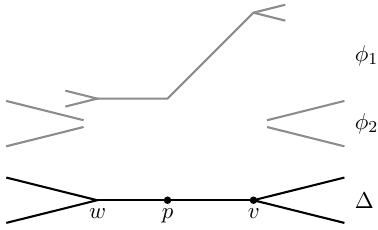}
\caption{PL functions $\phi_1$ and $\phi_2$ on an open cover of a graph $\Delta$, shown near a
single edge, used in the proof of Proposition~\ref{p:alg-equiv-curve} to show that the divisor $[p]$ is algebraically equivalent to
the divisor $[v]$.}\label{f:alg-equiv}
\end{figure}
First, suppose that $D$ is a divisor of degree $d$ on $\Delta$, and let $p$ be a
point of $D$ with multiplicity $m$. Suppose that $e$ is an edge containing $p$
and $v$ and $w$ are the endpoints of $e$. Define $U_1$ to be the open set
consisting of a small open neighborhood of $e$, and let $\phi_1$ be the PL
function on $U_1$ which is constant $0$ on $v$ and its incident edges, other
than $e$, has slope $m$ from $v$ to $p$, and is constant from $p$ to $w$ and on
the other edges incident to~$w$, as shown in Figure~\ref{f:alg-equiv}, in the
case $m =1$. Let $U_2$ be the complement of $e$, and
$\phi_2$ the constant $0$ function on $U_2$. Then, $\phi_1 \vert_{U_1 \cap U_2}
- \phi_2 \vert_{U_1 \cap U_2}$ is $0$ on the edges incident to $0$ and has a
constant, possibly non-zero value on the edges incident to $w$. Thus, these
functions define an algebraic equivalence between $D$ and $D-m[p] + m[v]$.

Note that $p$ may be the same as the endpoint $w$, and therefore repeating the
above construction means that $D$ is algebraically equivalent to $D - m[p] +
m[v']$
where $v'$ can be any vertex of $\Delta$. Repeating for all the points in $D$,
we get an algebraic equivalence between $D$ and $d v'$, where $d$ is the degree
of $D$. Therefore, all divisors of degree $d$ are algebraically equivalent to
each other.

Conversely, we want to show that algebraically equivalent divisors on~$\Delta$
have the same degree. More generally, we have the following proposition, which
applies to weak tropical complexes of arbitrary dimension.
\end{proof}

\begin{prop}\label{p:intersection-alg-equivalent}
If $D$ and $D'$ are algebraically equivalent divisors on a weak tropical
complex~$\Delta$, and $C$ is a curve, then $\deg D \cdot C = \deg D' \cdot C$.
\end{prop}

\begin{proof}
By the linearity of the intersection product, it is sufficient to show that the
degree of $D \cdot C$ is zero for $D$ an algebraically trivial divisor. Then,
there exist local defining equations for $D$ whose differences are locally
constant. We can subdivide $C$ into segments such that the defining equations of
$D$ are linear on the interior of each segment. The multiplicity of $D$ at a
point $p$ is, by definition, the sum of the outgoing slopes of these equations
along the segments containing~$p$. However, since the defining equations differ
by locally constant functions, the outgoing slopes at opposite ends of any
segment will cancel, giving that the total degree of $D \cdot C$ is $0$.
\end{proof}

While we expect that algebraically equivalence on the generic fiber $\X_\eta$ of
a degeneration will have a relationship with algebraic equivalence on the dual
complex~$\Delta$, we do not address that question in this paper. Instead, we
treat algebraic equivalence as a combinatorial definition, useful for
understanding the structure of the Picard group.

In particular, the remainder of this section will be used to show that the
N\'eron-Severi group is finitely generated, and the Jacobian is the quotient of
a real vector space by a discrete subgroup, proving
Theorem~\ref{t:picard-structure}. This is similar to the case of tropical
curves, where the Jacobian is known to be the quotient of $\RR^g$ by a
lattice~\cite{mikhalkin-zharkov}*{Thm.~6.2} (see \cite{baker-faber}*{Thm.~3.4}
for another proof), forming a topological torus of dimension equal to the first
Betti number of $\Delta$.

We begin by putting (\ref{eq:jac-pic-ns}) in the context of sheaf cohomology. We
define $\cA$ to be the sheaf of linear functions on~$\Delta$, and use $\RR$ to
denote the sheaf of locally constant functions on~$\Delta$, which is a subsheaf
of $\cA$. If we write $\cD$ for the quotient sheaf $\cA/\RR$, then we have a
long exact sequence on cohomology:
\begin{equation}\label{eq:exp}
\rightarrow
H^0(\Delta, \cD) \rightarrow H^1(\Delta, \RR) \rightarrow
H^1(\Delta, \cA) \rightarrow H^1(\Delta, \cD) \rightarrow
H^2(\Delta, \RR) \rightarrow
\end{equation}
This long exact sequence appeared in~\cite{mikhalkin-zharkov}*{Sec.~5.1},
where it was observed that $H^1(\Delta, \cA)$ can be identified with the Picard
group of~$\Delta$, and that (\ref{eq:exp}) has a striking similarity to the
exponential sequence on a complex algebraic variety. This exact
sequence~(\ref{eq:exp}) was subsequently used in~\cite{jell-rau-shaw} to prove
an analogue of the Lefschetz $(1,1)$-theorem for tropical manifolds.

\begin{prop}\label{p:picard-cohomology}
The Picard group of a weak tropical complex~$\Delta$ is isomorphic to
$H^1(\Delta, \cA)$. Moreover, the short exact at the middle term of
(\ref{eq:exp}) is isomorphic to (\ref{eq:jac-pic-ns}), meaning that
$\Jac(\Delta)$ is isomorphic to the image of $H^1(\Delta, \RR)$ in $H^1(\Delta,
\cA)$, and $\NS(\Delta)$ is isomorphic to the image of $H^1(\Delta, \cA)$ in
$H^1(\Delta, \cD)$.
\end{prop}

\begin{proof}
We let $\cP$ be the sheaf whose sections on a set $U$ consist of the piecewise
linear functions $\phi$ on $U$ which have only finitely many domains of
linearity, when restricted to a sufficiently small neighborhood of any point of
$U$. Then, $\cP$ is the sheaf associated to the presheaf of PL functions, which
were defined to have only finitely many domains of linearity. Note that, since
$\Delta$ is compact, a global section of $\cP$ is the same as a PL function.

We claim that the group of Cartier divisors on $\cP$ is isomorphic the group of
global sections of the quotient sheaf $\cP/\cA$. By definition any Cartier
divisor is defined by a collection of PL functions on an open cover of $\Delta$,
such that on the intersections, the PL functions define the same divisor, which
means their differences are linear by Proposition~\ref{p:divisor-pl}. On the
other hand, a global section of $\cP/\cA$ is given by a collection of sections
of $\cP$ on an open cover of $\Delta$, whose differences are sections of
$\cA$. Therefore, it is immediate that any Cartier divisor defines a global
section of $\cP/\cA$ and for the converse it suffices to show that the local
sections of $\cP$ each have only finitely many domains of linearity. This is
because the differences between global sections of are linear. Therefore, the
domains of linearity will agree between different local sections of $\cP$, so on the
compact topological space $\Delta$, there can only be finitely many domains of
linearity.

In addition, the principal divisors correspond to the elements of $\cP/\cA$
defined by a single global section of $\cP$.
Therefore, by the
long exact sequence in cohomology:
\begin{equation*}
H^0(\Delta, \cP) \rightarrow H^0(\Delta, \cP/\cA) \rightarrow H^1(\Delta, \cA)
\rightarrow H^1(\Delta, \cP),
\end{equation*}
the Picard group is the image of $H^0(\Delta, \cP/\cA)$ in $H^1(\Delta, \cA)$,
and to show that it is equal to $H^1(\Delta, \cA)$,
it will be sufficient to show that $H^1(\Delta, \cP)$ vanishes.

In fact, we will show that all higher sheaf cohomology of $\cP$ vanishes by
showing that $\cP$ is a soft sheaf~\cite{godement}*{Thm. II.4.4.3}. Recall that
$\cP$ is called a soft sheaf if for any closed set $Z \subset \Delta$, the map
$H^0(\Delta, \cP) \rightarrow H^0(Z, \cP)$ is surjective, where $H^0(Z, \cP)$ is
the direct limit of $H^0(U, \cP)$ as $U$ ranges over all open sets~$U$
containing~$Z$. We let $\phi$ be a function in $H^0(Z, \cP)$, which can be
represented by a PL function on some open set $U \supset Z$, which we also
denote~$\phi$. We can cover $Z$ by finitely many open cubes, each contained in
$U$. We dilate each of these cubes by a factor of $1-\epsilon$, such that their
union, which we denote $V$, still covers $Z$. In addition, the closure of $V$
will be a non-convex polyhedral set, which is also contained in $U$.

Since the closure $\overline V$ is compact, $\phi\vert_V$ is bounded, and we let $C$
be a constant less than the minimum of $\phi \vert_V$.
For any integer~$N$, we define the
function~$\phi_N$ on~$U$ by
\begin{equation*}
\phi_N(x) = \max\{C, \phi(x) - N d_1(\overline V, x)\},
\end{equation*}
where $d_1(\overline V,x)$ denotes the minimum distance, in the $L^1$~metric,
between $x$ and the closure of~$V$. It is clear that $\phi_N$ agrees with~$\phi$
on~$V$, and therefore, they have the same image in $H^0(\Delta, Z)$. Moreover,
$\phi_N$ is a PL function. For sufficiently large~$N$, we can ensure that
$\phi_N$ takes the value~$C$ at every point in a neighborhood of $\Delta
\setminus U$. Thus, we can extend $\phi_N$ by~$C$ to a function on all
of~$\Delta$, which completes the proof that $\cP$ is soft.

The assumptions on the local defining equations of an algebraically trivial
Cartier divisor~$D$ in Definition~\ref{d:alg-trivial-equiv} means that $D$
defines a global section in $H^0(\Delta, \cP/\RR)$. Using the long exact
sequence for the quotient $\cP/ \RR$, we get an element of $H^1(\Delta, \RR)$.
By the compatibility of long exact sequences, the divisor class of $D$ in
$H^1(\Delta, \cA)$ is in the image of $H^1(\Delta, \RR)$. Conversely, any
element of $H^1(\Delta, \RR)$ is in the image of $H^0(\Delta, \cP/\RR)$ by the
vanishing of $H^1(\Delta, \cP)$, and thus the classes of algebraically trivial
divisors in $H^1(\Delta, \cA)$ is exactly the image of $H^1(\Delta, \cA)$, which
proves the second statement of the proposition.
\end{proof}

\begin{prop}\label{p:differentials-discrete}
The image of $H^0(\Delta, \cD)$ in $H^1(\Delta, \RR)$ is discrete.
\end{prop}

\begin{proof}
It is sufficient to show that for any $\omega \in H^0(\Delta, \cD)$ and any
cycle $\gamma \in H_1(\Delta, \ZZ)$, the cap product $\omega \cdot \gamma \in
H_0(\Delta, \RR) \isom \RR$ is either $0$ or bounded away from $0$. The cycle
$\gamma$ can be represented by a formal sum of directed edges of $\Delta$, and
then $\omega \cdot \gamma$ is the sum of the slopes of $\cD$ along these edges.
Since sections of $\cD$ are required to have integral slopes, the pairing
$\omega \cdot \gamma$ will be an integer, which proves the proposition.
\end{proof}

By Propositions~\ref{p:picard-cohomology} and~\ref{p:differentials-discrete},
$\Jac(\Delta)$ is the
quotient of the real vector space $H^1(\Delta, \RR)$ by a discrete subgroup,
which therefore proves the second part of Theorem~\ref{t:picard-structure}. We
expect that in many cases of tropical complexes coming from degenerations,
$\Jac(\Delta)$ will be a topological torus, and even a tropical Abelian variety,
in the sense of~\cite{mikhalkin-zharkov}*{Sec.~5}. For a combinatorial result in
that direction, see~\cite{cartwright-surfaces}*{Thm.~1.2}.

\begin{figure}
\includegraphics{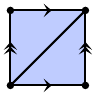}
\caption{The triangulation of the $2$-dimensional torus used in
Example~\ref{ex:torus}. A torus is formed by identifying the horizontal edges
with each other and the vertical edges with each other as indicated
by the arrow labeling.}\label{fig:torus}
\end{figure}

\begin{ex}\label{ex:torus}
As an example of using Proposition~\ref{p:picard-cohomology} to understand the
divisors on a tropical complex, we take $\Delta$ to be the triangulation of a
$2$-dimensional torus depicted in Figure~\ref{fig:torus} with $\alpha(v,e) = 1$
for all endpoints~$v$ of all edges~$e$. 
This example is the $2$-dimensional case of the theory of tropical Abelian
varieties discussed in~\cite{mikhalkin-zharkov}*{Section~5.1}.
The sheaf $\cD$ is isomorphic to the sheaf of locally constant
functions valued in~$\ZZ^2$, by taking a linear function to its
derivatives in the $x$ and~$y$ directions, and therefore $H^i(\Delta, \cD)$ is a
free Abelian group of rank $2 \binom{2}{i}$ for any $i$.

We claim that the long exact sequence~(\ref{eq:exp}) for $\Delta$ is:
\begin{equation*}
\rightarrow \ZZ^2 \rightarrow \RR^2 \rightarrow
(\RR/\ZZ)^2 \oplus \ZZ^3 \rightarrow
\ZZ^4 \rightarrow \RR \rightarrow
\end{equation*}
The summand~$(\RR/\ZZ)^2$ of the Picard group corresponds to algebraically
trivial divisors, each of which is linearly equivalent to a unique divisor of
the form $[\pi_1^{-1}(s)] - [\pi_1^{-1}(0)] + [\pi_2^{-1}(t)] -
[\pi_2^{-1}(0)]$, where $\pi_1$ and~$\pi_2$ are the two coordinate projections
from~$\Delta$ to the cycle of length~$1$ and $s$ and~$t$ are arbitrary points on
the $1$-cycle and $0$ is its vertex. The N\'eron-Severi group is~$\ZZ^3$, whose
generators can be taken to be the three edges in Figure~\ref{fig:torus}. Thus,
$\NS(\Delta)$ is a proper subgroup of $H^1(\Delta, \cD) \isom \ZZ^4$, and the map to
$H^2(\Delta, \RR) \isom \RR$ is non-trivial.
\end{ex}

We now prove the finite generation of $\NS(\Delta)$. For that, we use the
following lemma:

\begin{lem}\label{l:finite-generation}
Let $U$ be an open subset of any finite simplicial complex, consisting of a
union of interiors of the simplices. Suppose that $\cF$ is a sheaf on $U$, such
that, for each simplex $s$ whose interior $s^o$ is contained in $U$, the
restriction $\cF \vert_{s^o}$ is a constant sheaf with values in a finitely
generated group. Then $H^i(U, \cF)$ is a finitely generated group for each $i
\geq 0$.

Moreover, if $U$ consists of the union of interior of a single simplex~$s$ and
of the simplices containing $s$, then $H^i(U, \cF) = 0$ for $i > 0$.
\end{lem}

\begin{proof}
Let $d$ be the largest dimension of a simplex contained in the support of $\cF$,
and we will proceed by induction on $d$. When $d=0$, then $\cF$ is a direct sum
of skyscraper sheaves at the finitely many vertices contained in $U$, and so
$H^0(U, \cF)$ is finitely generated and $H^i(U, \cF) = 0$ for $i > 0$, which
proves both statements in the lemma.

Now suppose that $d$ is
positive, and let $U_{d-1}$ be the union of all simplices of dimension at most
$d-1$ in $U$. Define $\cG$ to be the push forward of the restriction of
$\cF$ to $U \setminus U_{d-1}$, i.e.\ the sheaf whose sections on $V \subset U$
are defined to be $H^0(V \setminus U_{d-1}, \cF)$. Since $\cG$ is constant on
the interior of each simplex, so is $\cG/\cF$, which is supported on $U_{d-1}$,
so it has finitely generated cohomology by the inductive hypothesis. In
addition, the cohomology of $\cG$ is the cohomology of the disjoint union of
interiors of simplices, which are contractible, and so $H^i(U, \cG) = 0$ for $i
> 0$. In addition, $H^0(U, \cG)$ is the direct sum of
the sections of finitely many simplices, so it too is finitely generated.
Finally, using the long exact sequence in cohomology, the cohomology of $\cF$ is
finitely generated, which completes the induction step and thus the proof of the
first statement.

Second, suppose that $U$ is the union of the interior of a single simplex~$s$
and of the simplices containing $s$. Then we again proceed by induction, with
the base case already established, and define $\cG$ as above. Then, for $i > 0$,
$H^i(U, \cG) = 0$, as already noted, and $H^i(U, \cG/\cF) = 0$ by the inductive
hypothesis, so it is sufficient to show that $H^0(U, \cG) \rightarrow H^0(U,
\cG/\cF)$ is surjective. Let $p$ be a point in the interior of $s$ and then
since any neighborhood of $p$ intersects every simplex of $U$, and both $\cG$
and $\cG/\cF$ are constant on the interiors of each simplex, we have isomorphism
between the global sections and the stalks at $p$, namely: $H^0(U, \cG) \isom
\cG_p$ and $H^0(U, \cG/\cF) \isom (\cG/\cF)_p$. However, surjectivity of the
homomorphism of stalks $\cG_p \rightarrow (\cG/\cF)_p$ follows from surjectivity
of the homomorphism of sheaves, and so $H^0(U, \cG) \rightarrow H^0(U, \cG/\cF)$
is surjective, proving the inductive step and thus the desired vanishing.
\end{proof}

\begin{proof}[Proof of Theorem~\ref{t:picard-structure}]
As noted above, Propositions~\ref{p:picard-cohomology}
and~\ref{p:differentials-discrete} shows that $\Jac(\Delta)$ is isomorphic to a
quotient of a real vector space by a discrete subgroup, so it just remains to
show that $\NS(\Delta)$ is finitely generated. For that, we claim that $\cD$ is
constant on the interiors of the simplices of $\Delta$, which will show that
$H^1(\Delta, \cD)$, as well as its subgroup $\NS(\Delta)$, are finitely
generated.

The set of linear functions on any connected open subset~$U$ of $\Delta$ is
determined by the set of facets and ridges that $U$ meets. For a sufficiently
small neighborhood of any open subset of the interior of a simplex~$s$ of
$\Delta$, these are just the facets, these are just the facets and ridges
containing $s$. Therefore, $\cA$ and $\cD$ are constant on the interior of~$s$.
In addition, the sections of $\cD$ are determined by the slopes on the facets,
which are integers, and thus the sections of $\cD$ are always finitely
generated. Therefore, $H^1(\Delta, \cD)$ is finitely generated by
Lemma~\ref{l:finite-generation}, and $\NS(\Delta)$ is a subgroup of $H^1(\Delta,
\cD)$ by Proposition~\ref{p:picard-cohomology}, so it is also finitely
generated.
\end{proof}

\begin{rmk}\label{r:computation}
The last sentence of Lemma~\ref{l:finite-generation} also gives a way of
effectively computing the Picard and N\'eron-Severi groups of a weak tropical
complex. In particular, for each simplex $s$, we can choose a small open
neighborhood~$U_s$ which only intersects $s$ and the simplices containing it,
and such that the intersections between any subset of the neighborhoods also
consists of a contractible open subset of a single simplex and of the simplices
containing it. Thus, although these intersections are not unions of interiors of
simplices of~$\Delta$, they are homeomorphic to a union of interiors of a
simplicial complex, and so we can apply Lemma~\ref{l:finite-generation} so that
$H^i(U_{s_1} \cap \cdots \cap U_{s_k}, \cA) = 0$ for any $i > 0$ and simplices
$s_1, \ldots, s_k$. Therefore, by~\cite{godement}*{Cor.\ to Thm.~II.5.4.1}, the
cohomology $H^1(\Delta, \cA) \isom \Pic(\Delta)$ can be computed using the \v
Cech cohomology of the cover $\{U_s\}$.
\end{rmk}

\end{document}